\RequirePackage[l2tabu, orthodox]{nag}
\documentclass[a4paper,10pt,reqno]{article}
\usepackage[utf8]{inputenc}
\usepackage[a4paper,
  inner=3.5cm, outer=3.5cm,
  top=3.5cm, bottom=2.5cm]{geometry}

\usepackage[english]{babel}

\usepackage{csquotes}
\usepackage{graphicx,tikz}
\usepackage{paralist}

\usepackage{amsmath,amssymb,amsthm}
\usepackage{mathtools,mathrsfs,dsfont}

\usepackage{booktabs,setspace,microtype,color}
\usepackage[hidelinks,pdfpagelabels]{hyperref}

\makeatletter
\newtheoremstyle{defindented}
  {3pt}
  {3pt}
  {\addtolength{\@totalleftmargin}{1em}
   \addtolength{\linewidth}{-2em}
   \parshape 1 1em \linewidth}
  {}
  {\bfseries}
  {.}
  {.5em}
  {}
\makeatother

\makeatletter
\newtheoremstyle{thmindented}
  {3pt}
  {3pt}
  {\addtolength{\@totalleftmargin}{1em} \itshape
   \addtolength{\linewidth}{-2em}
   \parshape 1 1em \linewidth}
  {}
  {\bfseries}
  {.}
  {.5em}
  {}
\makeatother

\makeatletter
\newtheoremstyle{remindented}
  {3pt}
  {3pt}
  {\addtolength{\@totalleftmargin}{1em}
   \addtolength{\linewidth}{-2em}
   \parshape 1 1em \linewidth}
  {}
  {\itshape}
  {.}
  {.5em}
  {}
\makeatother

\theoremstyle{thmindented}
\newtheorem{proposition}{Proposition}[section]
\newtheorem{theorem}[proposition]{Theorem}
\newtheorem{lemma}[proposition]{Lemma}
\newtheorem{corollary}[proposition]{Corollary}

\theoremstyle{defindented}
\newtheorem{definition}[proposition]{Definition}
\newtheorem{conditions}[proposition]{Assumptions}
\theoremstyle{remindented}
\newtheorem{remark}[proposition]{Remark}

\newcommand{\IB}{\mathbb{B}}
\newcommand{\IC}{\mathbb{C}}
\newcommand{\IL}{\mathbb{L}}
\newcommand{\IN}{\mathbb{N}}
\newcommand{\IQ}{\mathbb{Q}}
\newcommand{\IR}{\mathbb{R}}

\newcommand{\ga}{\alpha}
\newcommand{\gb}{\beta}
\newcommand{\gc}{\gamma}
\newcommand{\gC}{\Gamma}
\newcommand{\gd}{\delta}
\newcommand{\gD}{\Delta}
\newcommand{\gep}{\varepsilon}
\newcommand{\gvt}{\vartheta}
\newcommand{\gT}{\Theta}
\newcommand{\gi}{\iota}
\newcommand{\gl}{\lambda}
\newcommand{\gp}{\varphi}
\newcommand{\gs}{\sigma}
\newcommand{\gS}{\Sigma}
\newcommand{\go}{\omega}
\newcommand{\gO}{\Omega}

\newcommand{\cD}{\mathcal{D}}
\newcommand{\cE}{\mathcal{E}}
\newcommand{\cF}{\mathcal{F}}
\newcommand{\cH}{\mathcal{H}}
\newcommand{\cM}{\mathcal{M}}
\newcommand{\cS}{\mathcal{S}}

\newcommand{\scA}{\mathscr{A}}
\newcommand{\scB}{\mathscr{B}}
\newcommand{\scD}{\mathscr{D}}
\newcommand{\scE}{\mathscr{E}}
\newcommand{\scF}{\mathscr{F}}
\newcommand{\scG}{\mathscr{G}}

\newcommand{\scI}{\mathscr{I}}
\newcommand{\scS}{\mathscr{S}}
\newcommand{\scT}{\mathscr{T}}

\newcommand{\diam}{\operatorname{diam}}
\newcommand{\dimh}{\dim_{\mathcal{H}}}  
\newcommand{\dims}{\dim_{\mathcal{S}}}  
\newcommand{\dimw}{\dim_{\mathcal{W}}}  
\newcommand{\dimuw}{\overline{\dim_{\mathcal{W}}}} 
\newcommand{\SG}{\operatorname{SG}}     

\newcommand{\ssq}{\subseteq}
\newcommand{\Di}{\operatorname{d}\!}
\newcommand{\Prob}{\mathbf{P}}
\newcommand{\PTE}[1]{\mathbf{E}\left[{#1}\right]}
\newcommand{\PTEp}[2]{\mathbf{E}_{{#2}}\left[{#1}\right]}
\newcommand{\Lip}{\operatorname{Lip}}
\newcommand{\supp}{\operatorname{supp}}
\newcommand{\id}{\operatorname{id}}
\newcommand{\gr}{\operatorname{gr}}

\title{The Einstein Relation on Metric Measure Spaces}
\author{Burghart, Fabian\\
    \small Department of Mathematics and Computer Science\\ \small Eindhoven University of Technology, The Netherlands\\
    \small\texttt{f.burghart@tue.nl} 
    \and 
    Freiberg, Uta Renata\\
    \small Department of Mathematics\\ \small TU Chemnitz, Germany\\
    \small\texttt{uta.freiberg@mathematik.tu-chemnitz.de}
}
\date{February 2025}

\begin{document}

\maketitle

\begin{abstract}
 This note is based on F. Burghart's master thesis at Stuttgart University from July 2018, supervised by Prof. Freiberg. 
 
 We review the Einstein relation, which connects the Hausdorff, local walk and spectral dimensions on a space, in the abstract setting of a metric measure space equipped with a suitable operator. This requires some twists compared to the usual definitions from fractal geometry.  The main result establishes the invariance of the three involved notions of fractal dimension under bi-Lipschitz continuous isomorphisms between mm-spaces and explains, more generally, how the transport of the analytic and stochastic structure behind the Einstein relation works. While any homeomorphism suffices for this transport of structure, non-Lipschitz maps distort the Hausdorff and the local walk dimension in different ways. To illustrate this, we take a look at H\"older regular transformations and how they influence the local walk dimension and describe the Einstein relation on graphs of fractional Brownian motions. We conclude by giving a short list of further questions that may help building a general theory of the Einstein relation.
 
 
\end{abstract}

When regarding an open bounded domain $U$ in $\IR^n$, the Einstein relation is an equation expressing that the geometric behavior -- expressed by the asymptotic scaling of mass for small balls -- is nicely compatible with the analytic structure given by the Dirichlet-Laplace operator $\gD$ on $U$ -- expressed by the asymptotic behaviour of its eigenvalue counting function -- and with the asymptotic velocity of the stochastic process induced by $\gD$, namely Brownian motion. 
With the development of analytic and stochastic theory on (mainly self-similar) fractals, it was also discovered that the same relation holds on some fractals, most prominently on the Sierpi\'nski gasket $\SG$. 

The main goal of this article is to provide a general framework for the Einstein relation. To achieve this, we consider metric measure spaces $(X,d_X,\mu_X)$, where $(X,d_X)$ is a complete, separable, locally compact, and path-connected metric measure space not consisting of only a single point with an everywhere supported Radon measure $\mu_X$ on it. 

The purpose of Section~\ref{sec:Dims} is to formally introduce the Hausdorff dimension $\dimh$, the spectral dimension $\dims$ and the walk dimension $\dimw$. For $\dimh$, we give a short sketch of its definition and some of its properties, including Hutchinson's theorem on the Hausdorff dimension of self-similar sets in $\IR^n$. The spectral dimension is first motivated by Weyl's result on the asymptotic growth of the Dirichlet-Laplacian's eigenvalue counting function and then defined for operators $A$ on $L^2(X,\mu_X)$ that satisfy certain conditions. These conditions also ensure that there exists an essentially unique Hunt process (that is, a strong Markov process with certain continuity properties, see Subsubsection~\ref{sssec:DFtoMP} below) having $A$ as its infinitesimal generator. The outline of this theory relating processes, operator semigroups, generators and Dirichlet forms is then given in Subsection~\ref{ssec:dimw} before the (local) walk dimension and then the Einstein relation are defined.

In Section~\ref{sec:Examples}, we begin by examining the above mentioned classical case of a domain in Euclidean space with Dirichlet-Laplace operator, and continue by presenting the constructions of the standard Laplace operator on the Sierpi\'nski gasket $\SG$ as well as the construction of the Brownian motion on $\SG$. Both constructions rely heavily on the self-similarity and on the fact that $\SG$ can be approximated by a sequence of finite graphs, usually called prefractals. This also requires a different approach to the walk dimension than the local one from Subsection~\ref{ssec:dimw}, as the vertex set of a graph is always discrete. To see why this graph-theoretic walk-dimension can not be directly adapted to metric measure spaces, we present a counterexample in Subsection~\ref{ssec:realline}.

Section~\ref{sec:ERMMS} begins by defining two different types of morphisms between mm-spaces, namely contractions and Lipschitz-maps, both of which give rise to a notion of isomorphy, mm-isomorphy and Lipschitz-isomorphy, respectively. We then continue to investigate how we can transport the structure needed for the Einstein relation alongside maps 
\[
  \gp:(X,d_X,\mu_X)\to (Y,d_Y,\mu_Y)
\]
and prove that the Einstein relation is invariant under Lipschitz-isomorphisms. We proceed by looking at H\"older continuous transformations and manage to proof upper bounds for the walk dimension and apply this to Brownian motions running on the graphs of independent fractional Brownian motion, which generates a family of examples where the Einstein relation holds with a constant factor different from 1, see Theorem~\ref{thm:BHdimw}.  

The concluding discussion contains several open questions which may guide future development of a general theory of the Einstein relation as an invariant of metric measure spaces.

\section{Fractal Dimensions and the Einstein Relation}\label{sec:Dims}

In this introductory section, we wish to briefly expose the ingredients of the Einstein relation -- the Hausdorff dimension $\dimh$, the spectral dimension $\dims$, and the walk dimension $\dimw$ -- and state some of their properties. 

\subsection{Hausdorff measure and Hausdorff dimension}\label{ssec:Hausdorff}

Although the concepts of Hausdorff measure and dimension are well-known, we give the definitions in the interest of completeness. In what follows, let $(X,d)$ be a metric space.
\begin{definition}[Hausdorff outer measure]
  For fixed $s\geq0$, any subset $S\ssq X$ and any $\gd>0$, let 
  \[
    \cH^s_\gd(S)
      :=\inf \left\{\sum_{i\in I}(\diam U_i)^s:
            |I|\leq\aleph_0,S\ssq\bigcup_{i\in I}U_i\ssq X,\diam U_i\leq\gd\right\},
  \]
  i.e. the infimum is taken over all countable coverings of $S$ with diameter at most $\gd$. The $s$-dimensional Hausdorff outer measure of $S$ is now defined to be
  \begin{equation}\label{eq:DHM}
    \cH^s(S):=\lim_{\gd\searrow0}\cH^s_\gd(S).
  \end{equation}
\end{definition}
Observe that the limit in \eqref{eq:DHM} exists or equals $\infty$, since $\cH^s_\gd(S)$ is monotonically non-increasing in $\gd$ and bounded from below by 0. Furthermore, it can be shown that $\cH^s$ defines a metric outer measure on $X$, and thus restricts to a measure on a $\gs$-algebra containing the Borel $\gs$-algera $\scB(X)$ (cf. \cite[p.54ff]{mattila1999geometry}). By definition, the obtained measure then is the $s$-dimensional Hausdorff measure which we will denote by $\cH^s$ as well. Note that for $\cH^s$ to be a Radon measure, i.e. locally finite and inner regular, $\cH^s(X)<\infty$ is sufficient.

In the special case of $(X,d)$ being an Euclidean space, Hausdorff measures interpolate between the usual Lebesgue measures $\gl^n$: For $s=0$, we have simply $\cH^0(S)=\# S$, whereas for any integer $n>0$, it can be shown that there exists a constant $c_n>0$ depending only on $n$ such that $\cH^n=c_n\gl^n$, where the constant is the volume of the $n$-dimensional ball of unit diameter.

It can be seen by simple estimates that the map $s\mapsto \cH^s(S)$ for fixed $S\ssq X$ is monotonically non-increasing. More specifically, if $\cH^s(S)$ is finite for some $s$ then it vanishes for all $s'>s$, and conversely, if $\cH^s(S)>0$ then $\cH^{s'}(S)=\infty$ for all $s'<s$. Therefore, there exists precisely one real number $s$ where $\cH^\cdot(S)$ jumps from $\infty$ to $0$ (by possibly attaining any value of $[0,\infty]$ there). This motivates the following definition of Hausdorff dimension:
\begin{definition}
  The Hausdorff dimension $\dimh(S)$ of $S\ssq X$ is defined as
  \[
    \dimh(S):=\inf\{s\geq0:\cH^s(S)<\infty\}.
  \]
\end{definition}
Due to the above discussion, we have the following equalities:
\begin{align*}
   \dimh(S)&=\inf\{s\geq0:\cH^s(S)<\infty\}=\inf\{s\geq0:\cH^s(S)=0\}\\
           &=\sup\{s\geq0:\cH^s(S)=\infty\}=\sup\{s\geq0:\cH^s(S)>0\},
\end{align*}
providing some alternative characterisations of the Hausdorff dimension.

We further collect some important facts. To this end, let $S,S'$ and $S_1,S_2,...$ be subsets of $X$ as before. Then, the following properties hold (cf. \cite[p.32f]{falconer2007fractal} for a discussion in the Euclidean setting; however all arguments adapt to our more general situation without complication):
\begin{description}
  \item[Monotonicity.] If $S\ssq S'$ then $\dimh(S)\leq\dimh(S')$. 
  \item[Countable Stability.] For a sequence $(S_n)_{n\geq1}$, we have the equality
  \[ 
    \dimh\Big(\bigcup_{n\geq1}S_n\Big)=\sup_{n\geq1} \dimh(S_n).
  \]
  \item[Countable Sets.] If $|S|\leq\aleph_0$ then $\dimh(S)=0$.
  \item[H\"older continuous maps.] If $(X',d')$ is another metric space and $f:X\to X'$ is $\ga$-H\"older continuous for some $\ga\in(0,1]$ then $\dimh(f(S))\leq \ga^{-1}\dimh(S)$. In particular, the Hausdorff dimension is invariant under a bi-Lipschitz transformation (i.e. an invertible map $f$ with H\"older exponent $\ga=1$ for both $f$ and $f^{-1}$).
  \item[Euclidean Case.] If $(X,d)$ happens to be an Euclidean space (or more generally a continuously differentiable manifold) of dimension $n$ and $S$ is an open subset then $\dimh(S)=n$.
\end{description}

We conclude this section by discussing Hutchinson's theorem about the Hausdorff dimension of self-similar sets. For this, we recall that a map $F:X\to X$ on a metric space $(X,d)$ is a strict contraction if its Lipschitz constant satisfies
\begin{equation}\label{eq:Lip}
  \Lip_F:=\sup_{\stackrel{x,y\in X}{x\neq y}}\frac{d(F(x),F(y))}{d(x,y)}<1.
\end{equation}
If the stronger condition $d(F(x),F(y))=\Lip_F d(x,y)$ holds for all $x,y\in X$, we call $F$ a similitude with contraction factor $\Lip_F$.
\begin{theorem}[Hutchinson, \cite{hutchinson1981fractals}]\label{thm:hutchinson}
  Let $\scS=\{\cS_1,\dots,\cS_N\}$ be a finite set of strict contractions on the Euclidean space $\IR^n$. Then there exists a unique nonempty compact set denoted by $|\scS|$ which is invariant under $\scS$, i.e.
  \[
    |\scS|=\bigcup_{i=1}^N \cS_i(|\scS|).
  \]
  Furthermore, assume that $\scS$ satisfies the open set condition (OSC), meaning that there exists a nonempty bounded open set $O\ssq X$ with the properties $\cS_i(O)\ssq O$ and $\cS_i(O)\cap \cS_j(O)=\emptyset$ for all $i,j=1,...,N$ with $i\neq j$. Also assume that the maps $\cS_i$ are similitudes with contraction factor $r_i\in(0,1)$. Then, $\dimh(|\scS|)$ is the unique solution $s$ to the equation
  \[
    \sum_{i=1}^N r_i^s=1
  \]
  and we have $0<\cH^s(|\scS|)<\infty$. 
\end{theorem}
The proof is based on two different ideas: Existence and uniqueness can be shown in any complete metric space $(X,d)$ by invoking the Banach fixed point theorem for the map $X\supseteq A\mapsto\bigcup_{i=1}^N \cS_i(A)\ssq X$ acting on the space $(\cF(X),d_H)$ of all non-empty compact subsets of $X$ equipped with Hausdorff distance:
\[
  d_H(A,A'):=\inf\left\{\gep>0:A'\ssq B(A,\gep)\ \text{ and }\ 
         A\ssq B(A',\gep)\right\}\ \ \text{ for }\ A,A'\in\cF(X).
\]
Here, we denote by $B(A,\gep)$ the open $\gep$-neighbourhood of $A$. 
Indeed, $(\cF(X),d_H)$ is a complete metric space if $(X,d)$ is. The statement about the Hausdorff dimension and measure relies on a rather easy upper estimate on $\dimh$ and an application of the mass distribution principle to get $\cH^s(|\scS|)>0$: 
\begin{lemma}[Mass distribution principle, Frostmann, 
{\cite[Theorem 8.8]{mattila1999geometry}}]
  For a Borel set $X\ssq\IR^n$, we have $\cH^s(X)>0$ if and only if there exists a Borel measure $\mu$ on $X$ such that 
  $\mu(B(x,r))\leq r^s$ for $x\in\IR^n$ and $r>0$.
\end{lemma}
While uniqueness and existence of $|\scS|$ in Theorem~\ref{thm:hutchinson} are still ensured for maps on a complete metric space, the open set condition is not sufficient for statements about the Hausdorff dimension, see \cite{schief1996self} for further discussion.

\subsection{Weyl asymptotics and spectral dimension}\label{ssec:dims}

The idea of introducing spectral dimension is inspired by Weyl's law for the eigenvalues of the Dirichlet-Laplace operator which we will discuss here shortly before defining a larger class of operators that have similar spectral properties and are infinitesimal generators of Markov processes.

\subsubsection{The classical case}

Given a bounded open domain $U\ssq\IR^n$, consider the Laplace operator $\gD$ acting on functions $u:U\to\IR$ which satisfy the Dirichlet boundary condition $u\equiv0$ on $\partial U$. Then, the spectrum of $-\gD$ consists of non-negative eigenvalues with a single accumulation point at $\infty$. Hence we can order them in a non-decreasing way, counting the geometric multiplicities, as
\begin{equation}\label{eq:EVs}
  0\leq\gl_1\leq\gl_2\leq...\leq\gl_n\leq...\ \text{ with }\ \gl_n\nearrow\infty. 
\end{equation}
In this setting, it makes sense to define the eigenvalue counting function via 
\begin{equation}\label{eq:DECF}
  N_{-\gD}(x):=\max\{n\in\IN:\gl_n\leq x\},\ x\in\IR_{\geq0}.
\end{equation}
Weyl's law now states that there is the asymptotic equivalence\footnote{We adopt the notation $f\sim g$ for the equivalence relation given by $\lim\frac{f}{g}=1$.}
\begin{equation}\label{eq:WL}
  N_{-\gD}(x)\sim C_n\cH^n(U)x^{n/2},\quad x\nearrow\infty,
\end{equation}
where the constant $C_n$ is independent of the domain $U$ (see \cite{Weyl1911} and \cite{Weyl1912} for the original publications). Motivated by \eqref{eq:WL}, we define the spectral dimension of $-\gD$ on $U$ by
\begin{equation}\label{eq:dims}
  \dims(U,-\gD):=\lim_{x\to\infty}\frac{\log N_{-\gD}(x)}{\log x}
\end{equation}
which yields $n/2$ in the situation examined by Weyl's law. Note that the usual definition of $\dims$ differs by a factor of 2 (cf. \cite{kigami1993weyl},\cite{hambly_kigami_kumagai_2002}) so that $\dims(U,\gD)$ normally coincides with $\dimh(U)=n$. However, this comes at the cost of an additional factor in the Einstein relation. Moreover, it can be argued that the spectral dimension is rather a property of the operator $-\gD$ than of the underlying space $U$; accordingly, it would perhaps be more accurate to call the quantity in \eqref{eq:dims} the spectral exponent of $-\gD$. Nonetheless, we stick to the name spectral dimension but take the liberty to deviate from the established convention in this minor aspect.

\subsubsection{The general case}

How can we generalise the concepts just introduced to sets which are not bounded open subsets of $\IR^n$? For this purpose, suppose we are given a metric measure space 
$(X,d,\mu)$, where $(X,d)$ is a locally compact separable metric space and $\mu$ is a Radon measure on $X$.

Of course, it is possible to define an eigenvalue counting function in the same fashion as above for any operator $A$ whose set of eigenvalues possesses only one accumulation point at $+\infty$. However, as we will explain in the next section, we also wish to associate a reasonably well-behaved Markov process with state space $X$ to $A$. This requires several additional assumptions that will be motivated in this and the next section. More precisely, we choose to impose the following conditions on $A$:
\begin{conditions}\label{cond:A}
For an operator $A:L^2(X,\mu)\supseteq\scD(A)\to L^2(X,\mu)$, we assume that the following statements hold:
\begin{description}
  \item[Self-adjointness.] $A$ is a densely defined, self-adjoint (and therefore closed) operator on the Hilbert space $L^2(X,\mu)$.
  \item[Eigenvalues.] The spectrum is contained in $\IR_{\geq0}$, nonempty, and has a single accumulation point at $+\infty$. Hence the set of eigenvalues can be enumerated as in \eqref{eq:EVs}.
  \item[Regularity (of the corresponding Dirichlet form).] The set\footnote{We denote by $C_c(X)$ the space of all compactly supported continuous functions on $X$} $\scD(\sqrt{A})\cap C_c(X)$ is dense in $C_c(X)$ with respect to the supremum norm and is dense in the domain $\scD(\sqrt{A})$ endowed with the graph norm 
  \[
    \|f\|_{\sqrt{A}}:=\|f\|_{L^2(X,\mu)}+\|\sqrt{A}f\|_{L^2(X,\mu)},\ \ f\in\scD(\sqrt{A}).
  \]
  \item[Dissipativeness.] $-A$ is dissipative. In other words, for all $f\in\scD(A)$ and all $\gl>0$, we have $\|(\gl+A)f\|_{L^2(X,\mu)}\geq\gl\|f\|_{L^2(X,\mu)}$.
\end{description}
\end{conditions}
Here, $\sqrt{A}$ denotes an operator satisfying $\sqrt{A}\circ\sqrt{A}=A$ that can be defined via a standard spectral-theoretic construction. We will not go into greater detail here and refer to \cite{birman2012spectral} instead.

The first of these assumptions guarantees that $A$ is a closed operator, whereas the second ensures that $\gl+A$ is surjective for at least one $\gl>0$. Thus, the Hille-Yosida theorem states that there is a strongly continuous semigroup of contractive linear operators $T_t$ on $H$ such that $-A$ is its infinitesimal generator. That is to say:
\begin{definition}\label{def:semigroup}
  A strongly continuous semigroup $(T_t)_{t\geq0}$ on a Hilbert space $H$ is a monoid homomorphism $t\mapsto T_t$ from $(\IR_{\geq0},+)$ to the space of bounded linear operators $(\IB(H),\cdot)$ on $H$ (equipped with composition) satisfying for all $f\in H$ the additional property
  \[
    \lim_{t\searrow0} \|T_tf-f\|=0.
  \]
  The infinitesimal generator $(-A,\scD(A))$ of $(T_t)_{t\geq0}$ is defined via
  \[
    (-A)f=\lim_{t\searrow0}\frac{1}{t}(T_tf-f),\ f\in\scD(A),
  \]
  where $\scD(A)$ is the set of elements in $H$ for which this limit exists.
\end{definition}
\begin{theorem}[Hille-Yosida,\cite{ma2012introduction}]\label{thm:HY}
  An operator $(-A,\scD(A))$ is the generator of a strongly continuous semigroup $(T_t)_{t\geq0}$ with $\|T_t\|\leq1$ for all $t\geq0$ if and only if $-A$ is a densely defined, closed, dissipative operator such that for some $\gl>0$, the map $\gl+A$ is surjective. 
\end{theorem}
It can be shown that there is a one-to-one correspondence between contractive semigroups and operators that satisfy the Hille-Yosida theorem, that is, the semigroup in the above theorem is uniquely determined by $A$. 

\begin{remark}
  In the above list of assumptions, dissipativeness is redundant as we regard operators on Hilbert spaces. In this setting, $-A$ is dissipative if $A$ is a positive operator since
  \[
    \|(\gl+A)f\|^2-\gl^2\|f\|^2=\|Af\|^2+2\gl\left<f,Af\right>\geq0
  \]
\end{remark}

Having discussed the motivation for the assumptions \ref{cond:A}, we now proceed to adapt the definitions made in \eqref{eq:DECF} and \eqref{eq:dims} in a rather straightforward way:
\begin{definition}
  Given an operator $(A,\scD(A))$ on $L^2(X,\mu)$ satisfying the assumptions \ref{cond:A}, its eigenvalue counting function is defined by
  \begin{equation}\label{eq:defXVCF}
    N_A(x):=\max\{n\in\IN:\gl_n\leq x\},\ x\in\IR_{\geq0},
  \end{equation}
  and, if the limit exists, the spectral dimension of $A$ by
  \begin{equation}\label{eq:defdims}
    \dims(X,A):=\lim_{x\to\infty}\frac{\log N_A(x)}{\log x}.
  \end{equation}
\end{definition}

\subsection{Markov processes and walk dimension}\label{ssec:dimw}

In this subsection, we review the connection between Markov processes, semigroups, Dirichlet forms, and their infinitesimal generators.

\subsubsection{From Dirichlet forms to Markov processes}\label{sssec:DFtoMP}

The theory presented here is mostly taken from \cite{fukushima2011dirichlet} and \cite[chapter 4]{ma2012introduction}. We work on $L^2(X,\mu)$ where $\mu$ is a $\gs$-finite Borel-measure on $X$. 
\begin{definition}\label{defin:DF}
  A map $\cE:\scD(\cE)\times\scD(\cE)\to\IR$ is a Dirichlet form if it satisfies the following conditions:
  \begin{enumerate}
  \item The domain $\scD(\cE)\ssq L^2(X,\mu)$ of $\cE$ is a dense linear subspace.
  \item $\cE$ is a symmetric, non-negative definite bilinear form.
  \item This form is closed, that is, the inner product space $(\scD(\cE),\cE_\ga)$ equipped with the scalar product
  \[
    \cE_\ga(u,v):=\cE(u,v)+\ga\left<u,v\right>\ \text{ for }\ u,v\in\scD(\cE_\ga)=\scD(\cE),\ \ga>0,
  \]
  is complete (and thus itself a Hilbert space).
  \item $\cE$ is a Markovian form, i.e. for all $u\in\scD(\cE)$, $v:=(0\vee u)\wedge 1\in\scD(\cE)$ and we have $\cE[v]\leq\cE[u]$ for the quadratic form of $\cE$. 
  \end{enumerate}
\end{definition}
We remark that the choice of $\ga>0$ is irrelevant for the completeness of $(\scD(\cE),\cE_\ga)$ since all induced norms are equivalent to each other. 
\begin{definition}
  A Dirichlet form $(\cE,\scD(\cE))$ on $L^2(X,\mu)$ is said to be
  \begin{enumerate}
    \item regular if it possesses a core, that is, the space $\scD(\cE)\cap C_c(X)$ is simultaneously dense in $\scD(\cE)$ with respect to the $\cE_1$-norm and in $C_c(X)$ with respect to the uniform norm. 
    \item local if $\cE(u,v)=0$ whenever $u,v\in\scD(\cE)$ have disjoint compact support.
    \item strongly local if $\cE(u,v)=0$ whenever $u,v\in\scD(\cE)$ have compact support and $v$ is constant on a neighbourhood of $\supp(u)$.
  \end{enumerate}
  If additionally $\mu(X)<\infty$, we say that $\cE$ is
  \begin{enumerate}
    \setcounter{enumi}{3}
    \item conservative if $1\in\scD(\cE)$ and $\cE[1]=0$.
    \item irreducible if it is conservative and $\cE[f]=0$ implies that $f$ is constant. 
  \end{enumerate}
\end{definition}
We can uniquely attach a positive semidefinite operator $A$ to a Dirichlet form (and vice versa) via the relation 
\begin{equation}\label{eq:DFtoOp}
  \cE(u,v)=\left<Au,v\right>,\  u\in\scD(A),v\in\scD(\cE).
\end{equation}
In particular, if $A$ meets the requirements of \ref{cond:A}, we not only have precisely one strongly continuous contraction semigroup on $L^2(X,\mu)$ as explained by Theorem~\ref{thm:HY}, but also a unique Dirichlet form thanks to \eqref{eq:DFtoOp}. In similar style, we would also like to attach a unique Markov process to $A$; or, equivalently, to the semigroup or the Dirichlet form. 

To define a suitable stochastic process with values in $X$, we first adjoin a cemetery state $*$ in such a way that if $X$ is non-compact, $X_*:=X\uplus\{*\}$ is the one-point compactification of $X$, whereas $*$ is supposed to be an isolated point if $X$ is compact. Let $M=\left(\gO,\scA,(M_t)_{t\geq0},(\Prob_x)_{x\in X_*}\right)$ be a stochastic process on a measurable space $(\gO,\scA)$ with values in $X_*$, where we adapt the notation that $\Prob_x[M_0=x]=1$ for all 
$x\in X_*$ and $\Prob_*[M_t=*]=1$ for all $t\geq0$ (where the latter means that the cemetery state $*$ is absorbing). Note that $M$ induces a filtration $\scF=(\scF_t)_{t\geq0}$ on $\scA$ by 
\[
  \scF_t=\bigcap_{\Prob\in\cM^+_1(\gO,\scA)} \big(\gs\{M_s:0\leq s\leq t\}\big)^\Prob.
\]
Here, $\cM^+_1$ denotes the set of all probability measures on $(\gO,\scA)$, $\gs\{\cdot\}$ denotes the the $\gs$-algebra generated by $\{\cdot\}$ and $\scB^\Prob$ denotes the completion of a $\gs$-algebra $\scB$ with respect to the measure $\Prob$. Henceforth, we will only consider stochastic processes $M$ that satisfy the strong Markov property with respect to $\scF$ and are time-homogeneous. Such $M$ is called a Hunt process if it additionally has right-continuous trajectories with left-limits and is quasi-left-continuous, i.e.  any sequence $\tau_n\nearrow\tau$ of $\scF$-stopping times satisfies
\[
  \Prob_\ga\left[\lim_{n\to\infty} M_{\tau_n}=M_\tau,\tau<\infty\right]=\Prob_\ga[\tau<\infty]
\]
for any initial distribution $\ga$. We can now translate Markov processes to contractive semigroups by setting
\begin{equation}\label{eq:MPtoOp}
  (T_tf)(x):=\PTEp{f(M_t)}{x},\ t\geq0.
\end{equation}
The other direction is more involved, and the process attached to a Dirichlet form is generally non-unique. We have, however, (cf. \cite[Theorems 7.2.1 and 7.2.2]{fukushima2011dirichlet})
\begin{theorem}\label{thm:fukushima}
  Let $\cE$ be a regular Dirichlet form on $L^2(X,\mu)$. Then, there exists a Hunt process $M$ on $(X,d)$ such that the operators $T_t, t\geq0$, from \eqref{eq:MPtoOp} are symmetric and $\cE$ is the Dirichlet form belonging to this semigroup.
  
  Moreover, if $\cE$ is local, $M$ is a diffusion process.
\end{theorem}
As hinted above, those processes are not unique: One can modify $M$ to $\tilde M$ by killing the process on a polar set and obtain the same semigroup for both. See Section~7.2.2. in \cite{fukushima2011dirichlet} for further discussion.

\subsubsection{Local walk dimension and Einstein relation}\label{sssec:dimwER}

The walk dimension is meant to quantify how fast a given Markov process on $M$ moves away from its starting point $x$. This is best expressed in terms of the stopping time $\tau(r)=\tau(x,r):=\tau(B(x,r))$, which is supposed to be the first exit time of the ball $B(x,r)=\{y\in X:d(x,y)<r\}$ (we generally suppress the dependency on the starting point). Note that for a Hunt process, this is indeed an $\scF$-stopping time by the right-continuity of the process in question and by \cite[Lemma 7.6]{kallenberg2002foundations}.

For the next definition to make sense, we need to impose some additional assumption on the metric space $(X,d)$. We choose to demand that $X$ is path connected and does not consist of a single point, but will also discuss the case where it is the vertex set of a graph in the next section. 
\begin{definition}\label{def:dimw}(Cf. \cite{hambly_kigami_kumagai_2002})
  We define the quantity
  \[
    \dimw(X,M;x)=\lim_{r\searrow0}\frac{\log \PTEp{\tau(r)}{x}}{\log r}
  \]
  and call it the (local) walk dimension of $(X,d)$ at $x\in X$ with respect to the Markov process $(M_t)_{t\geq0}$. If $\dimw(X,M;x)$ is $\mu$--a.e. constant on $X$, we shorten our notation to $\dimw(X,M)$.
\end{definition}
Of course, this definition makes sense for almost any stochastic process, so whenever we are interested in the walk dimension alone, we do not need to assume that the process $M$ is a Hunt process.

We are now finally able to state the Einstein relation: 
\begin{definition}\label{defin:ER}
  Let $(X,d,\mu)$ be a locally compact separable metric measure space and let $(A,\scD(A))$ be an operator on $L^2(X,\mu)$ satisfying assumptions \ref{cond:A}. Suppose $M=\big((M_t)_{t\geq0},(\Prob_x)_{x\in X_*}\big)$ is a Markov process associated to $A$ via the Dirichlet form $\cE(\cdot,\cdot)=\left<\sqrt{A}\,\cdot,\sqrt{A}\,\cdot\right>$. We then say that the Einstein relation with constant $c$ holds on $X$ with respect to $A$ if
  \begin{equation}\label{eq:ER}
    \dimh(X)=c\dims(X,A)\dimw(X,M).
  \end{equation}
  We omit mentioning the constant if $c=1$.
\end{definition}
Of course, we are mainly interested in the case $c=1$ since this means that geometry, analysis and stochastic on the given mm-space ``fit well together''. Nonetheless, the invariance properties derived in Section~\ref{sec:ERMMS} hold regardless of the concrete value of this constant.

\subsection{Related works}\label{ssec:Related}

The Einstein relation exists in several different version, each adapted to its setting. 

On graphs, where the edges are interpreted as unit resistors in an electric network, it can be formulated in terms of the
volume growth rate of balls, the mean exit time of a random walk from a ball centered at its starting point and the growth rate of the resistance of an annulus. Then, the Einstein relation holds if all three quantities are well-defined and the first one is the sum of the other two. See \cite{telcs2006art} for the exact definitions and the theory in this setting. 

For fractals, the Einstein relation is usually considered on post-critically finite, self-similar fractals. To elaborate, note that in the setting and notation of Theorem~\ref{thm:hutchinson}, $|\scS|$ is obtained from $\scS=\{\cS_1,\dots,\cS_N\}$ in such a way that points $x\in|\scS|$ can be encoded by their ``address'' $w\in\{1,\dots,N\}^\IN$ through a quotient map $\pi:\{1,\dots,N\}^\IN \twoheadrightarrow |\scS|: w\mapsto x$ with $x$ being the unique point in $(\cS_{w_1} \circ \cS_{w_2} \circ \dots)(|\scS|)$. Writing $C:= \pi^{-1}\left(\bigcup_{i,j=1,\dots,N;\ i\neq j} \cS_i(|\scS|) \cap \cS_j(|\scS|)\right)$ for the set of addresses belonging to overlap points between the $\cS_i(|\scS|)$ and $\gs:\{1,\dots,N\}^\IN \to \{1,\dots,N\}^\IN:w_1w_2w_3\dots \mapsto w_2w_3\dots$ for the usual shift map, a self-similar fractal is post-critically finite if the set $\bigcup_{n\geq 1} \gs^n(C)$ is finite, see \cite[Definition~1.3.13]{kigami2001analysis}. 
In this setting, the usual Einstein relation is very similar to the one introduced in Definition~\ref{defin:ER}. Since these spaces admit a meaningful approximation by prefactals, i.e. finite graphs, one can define the walk dimension globally as 
\begin{equation}\label{eq:defdimwrong}
  \lim_{R\nearrow\infty}\frac{\log \PTEp{\tau(R)}{x}}{\log R},
\end{equation}
which is essentially the same definition as for graphs (under certain conditions, this limit is independent of the starting point in graphs with infinite diameter). Using the self-similarity of the fractal, one can hope to avoid taking the limit in \eqref{eq:defdimwrong}, see \cite{freiberg2012einstein}. We return to this setting in Subsection~\ref{ssec:SG} below. In another context, the walk dimension occurs as an exponent in heat kernels, see among others \cite{barlow1998diffusions} and \cite{grigoryan2021analysis}.

Finally, we mention that \cite{hambly_kigami_kumagai_2002} considered a completely localised version of the Einstein relation for a multifractal formalism. This variant relied on the local walk dimension of Definition~\ref{def:dimw} yet also featured a local geometric dimension based on a given measure and a local spectral dimension, defined via estimates to the transition kernel of the operator semigroup.

\section{Examples and Non-examples}\label{sec:Examples}

In this section, we will discuss the necessity of some of the restrictive assumptions made previously and explore the Einstein relation by examining some examples and -- by doing so -- will motivate some of the more general results of the next section.

\subsection{Euclidean Space}\label{ssec:Euclidean}

We start by examining the classical setting of paragraph 1.2.1 in greater detail: Let once again $U\ssq\IR^n$ be an open, bounded, non-empty domain, equipped with Euclidean metric and $n$-dimensional Lebesgue-measure $\gl^n$. Trivially, $\dimh(U)=n$. 

For the Dirichlet-Laplace operator as introduced earlier, we obtain 
$\dims(U,-\gD)=\frac{n}{2}$ due to \eqref{eq:DECF}-\eqref{eq:dims}. Simultaneously, it is well-known that $\frac{1}{2}\gD$ is the generator of the $n$-dimensional Brownian motion $B_t$, which can easily be seen as follows:

By \eqref{eq:MPtoOp}, the semigroup $(T_t)_{t\geq0}$ induced by $B_t$ reads
\[
  T_tf(x)
   =\PTEp{f(B_t)}{x}
   =\frac{1}{(2\pi t)^{n/2}}\int_{\IR^n} 
     \exp\left(-\frac{|x-y|^2}{2t}\right)f(y)\,\Di y,
\]
for $f\in L^2(\IR^n,\gl^n)\cap C_c(E)$.
Comparing this expression to the well-known convolution formula (see \cite[p.47]{evans2010partial})
\[
  u(x,t)=\frac{1}{(4\pi t)^{n/2}}\int_{\IR^n} \exp\left(-\frac{|x-y|^2}{4t}\right)f(y)\,\Di y,\ x\in\IR^n, t\geq0,
\]
for the solution of the heat equation $\gD u=\partial_t u$ on $\IR^n$ with initial value $u(x,0)=f(x)$, we can quickly derive that $\gD T_tf(x)=2\partial_t T_tf(x)$. Since imposing the Dirichlet boundary conditions on $\gD$ corresponds to killing $B_t$ at the boundary of $U$, we thus obtain by Definition~\ref{def:semigroup} the generator
\[
  Af=\frac{1}{2}\gD f,
\]
extended to its maximal domain, the Sobolev space $H_0^1(U,\mu)$. We conclude that the Markov process associated with $\frac{1}{2}\gD$ is $B_t$.

It remains to determine the walk dimension of the $n$-dimensional Brownian motion. This can be done in several ways, for example by appealing to Brownian scaling or by invoking Dynkin's formula after a standard truncation argument: By applying \cite[Lemma 19.21]{kallenberg2002foundations} to the function $u_x(y)=\frac{1}{2n}|y-x|^2$, we get
\[
  \PTEp{\tau(r)}{x}
  =\PTEp{\int_0^{\tau(r)}\gD u_x(2B_s)\,\Di s}{x}
  =\PTEp{u_x\left(B_{\tau(r)}\right)-u_x(0)}{x}
  =\frac{r^2}{n}.
\]
Therefore, by Definition~\ref{def:dimw}, we obtain $\dimw(U,B_t)=2$ which implies together with the results obtained previously that the Einstein relation (with constant 1) holds on $U$ with respect to $\gD/2$.

\subsection{The Sierpi\'nski Gasket}\label{ssec:SG}

\begin{figure}[ht]
\centering
\includegraphics[scale=7]{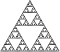}
\caption{An approximation of $\SG$}\label{fig:sg}
\end{figure}
The Sierpi\'nski Gasket is a simple example of a fractal obtained by means of an iterated function system and can be described according to Theorem~\ref{thm:hutchinson} as the unique non-empty compact set $\SG\ssq\IR^2$ which is invariant under the three similitudes 
\[
  \cS_1(x,y)=\left(\frac{x}{2},\frac{y}{2}\right),\ 
  \cS_2(x,y)=\left(\frac{x+1}{2},\frac{y}{2}\right),\ 
  \cS_3(x,y)=\left(\frac{2x+1}{4},\frac{2y+\sqrt{3}}{4}\right),
\]
see figure \ref{fig:sg}. Since $\{\cS_1,\cS_2,\cS_3\}$ satisfies the (OSC), e.g. by taking the open equilateral triangle with corners $(0,0), (0,1)$ and $(1/2,\sqrt{3}/2)$, we obtain both 
\begin{equation}\label{eq:SGdimh}
  s=\dimh(\SG)=\frac{\ln 3}{\ln 2}
\end{equation}
and $\cH^s(\SG)\in(0,\infty)$ by a second appeal to Hutchinson's theorem. 

We will use the remainder of this section to establish the validity of the Einstein relation on $\SG$ with respect to the standard Laplace operator on $\SG$, which can be obtained in two different ways. In order to describe these constructions, we first need to fix some notation. 

Let $\scS=\{\cS_1,\cS_2,\cS_3\}$ as in section 1.1, set $\gS=\{1,2,3\}$ and denote by 
\[
  \gS^*:=\{\gep\}\cup\bigcup_{n\geq 1}\gS^n
\]
the free monoid consisting of all finite words over the alphabet $\gS$, where the monoid operation is given by concatenation and $\gep$ is supposed to represent the empty word. Given a word of length $l$, say $w=w_1...w_l\in\gS^l$, we define the shorthand notation $S_w:=S_{w_1}\circ\dots\circ S_{w_l}$. By abuse of notation, we will simply write $w(x)$ instead of $S_w(x)$. By an $n$-cell we understand the set $w(\SG)\ssq\SG$, where $w$ is a word of length $n$. Note that two different cells are either disjoint, or intersect in a single point which we will then call conjunction point, or one of them is completely contained in the other one.

\begin{figure}[ht]
\centering
\includegraphics[width=0.95\textwidth]{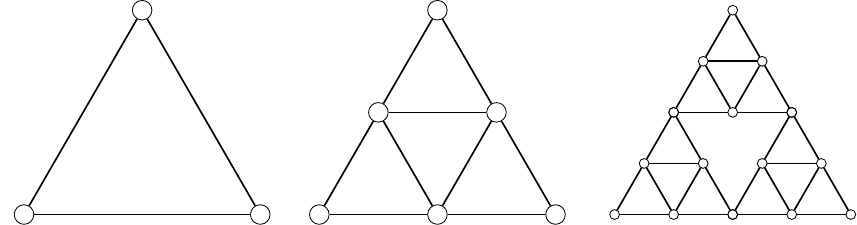}
\caption{The graphs $G_0,G_1,G_2$}\label{fig:gn}
\end{figure}

It is possible to approximate $\SG$ by a sequence of graphs $G_n$. Those graphs can be thought of as planar graphs with a triangle for each $n$-cell, where the vertices are the conjunction points between them, compare Figure~\ref{fig:gn}. More precisely, let $G_n$ be the graph embedded in $\IR^2$ with vertex set $V_n$ inductively defined by
\begin{align*}
  V_0&:=
  \left\{(0,0),(0,1),\left(\frac{1}{2},\frac{\sqrt{3}}{2}\right)\right\}\\
  V_{n+1}&:=S_1(V_n)\cup S_2(V_n)\cup S_3(V_n),\ \ n\geq0.
\end{align*}
Additionally, set 
\[
  V^*:=\bigcup_{n\geq0} V_n.
\]
Note that $V_0\ssq V_1\ssq...\ssq V_n\ssq ...$ and that
\[
  V_n=\bigcup_{w\in\gS^n} w(V_0).
\]
It also follows from the proof of Hutchinson's theorem that 
$d_H(V_n,\SG)\to 0$ as $n\to\infty$, so indeed, the graphs $G_n$ approximate $\SG$. 

In $G_n$, connect two vertices $x,y\in V_n$ by a straight edge $xy$ if and only if they belong to the same $n$-cell, cf. Figure~\ref{fig:gn}, in which case we will call them neighbours and write $x\sim_n y$. By $E_n$ we mean the set of all edges in $G_n$.

\subsubsection[Approximation by Dirichlet forms]{Approximation by Dirichlet forms\protect\footnote{The material in this section is an overview of the construction given in \cite[chapter I]{strichartz2006differential}}}

This analytic approach works by establishing so-called energy forms on graphs $G_n$ -- these are graph-theoretic discretisations of the Dirichlet form attached to the Laplace operator. For $n\in\IN$, define a bilinear form $\tilde\cE_n$ on $L^2(V_n)\cong \IR^{V_n}$ by 
\[
  \tilde\cE_n(f,g):=\sum_{xy\in E_n}(f(x)-f(y))(g(x)-g(y)), 
  \ \ \ f,g\in L^2(V_n).
\]
It is easy to check that this defines a local Dirichlet form on $L^2(V_n)$. As we will see, these bilinear forms are compatible with each other after a suitable renormalisation. Starting from $\tilde\cE_0$, consider a function 
$u\in\IR^3\cong L^2(V_0)$. Under all possible extensions 
$\tilde u$ to $\IR^6\cong L^2(V_1)$, which function is the harmonic extension of $u$, i.e. minimises $\tilde\cE_1[\tilde u]$? 

\begin{figure}[ht]
\centering
\includegraphics[scale=1]{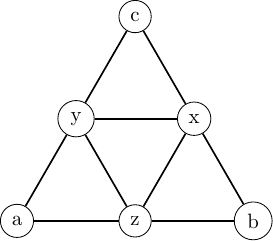}
\caption{The values of $u$ on $V_1$}\label{fig:g1}
\end{figure}

To answer this question, we label each vertex in $V_0$ by its value under $u$, say $a,b,c$ and each vertex in $V_1\setminus V_0$ by its value under $\tilde u$, say $x,y,z$ as in figure \ref{fig:g1}. Since we assume $u$ to be a fixed given function, the values $a,b,c$ are fixed, and we obtain 
\begin{multline}\label{eq:E1}
  \tilde\cE_1[\tilde u]=\left[(x-b)^2+(x-c)^2+(y-a)^2+(y-c)^2+(z-a)^2+(z-b)^2\right.\\
  \left.+(x-y)^2+(y-z)^2+(z-x)^2\right]
\end{multline}
Finding the minimising values for $x,y,z$ now becomes an exercise in multivariate calculus. Setting the partial derivatives equal to 0 yields the following system of linear equations:
\begin{align*}
  \begin{array}{c}
    8x=2(b+c+y+z)\\
    8y=2(a+c+x+z)\\
    8z=2(a+b+x+y)\\
  \end{array}\Longleftrightarrow
  \begin{array}{c}
    5x=a+2b+2c\\
    5y=2a+b+2c\\
    5z=2a+2b+c\\
  \end{array}
\end{align*}
Plugging this solution as $\tilde u$ in \eqref{eq:E1} to evaluate $\tilde\cE_1[\tilde u]$ gives
\begin{align*}
  \tilde\cE_1[\tilde u]
  &=\frac{1}{25}\left[(a-3b+2c)^2+(a+2b-3c)^2+(-3a+b+2c)^2
    +(2a+b-3c)^2\right.\\
  &\qquad \left.+(-3a+2b+c)^2+(2a-3b+c)^2
    +(b-a)^2+(c-b)^2+(a-c)^2\right]\\
  &=\frac{30}{25}\left[a^2+b^2+c^2-ab-ac-bc\right]\\
  &=\frac{3}{5}\tilde\cE_0[u]
\end{align*}
Using the fact that the vertices in $V_{n+1}\setminus V_n$ are in $G_n$ only adjacent to vertices in $V_n$ and that 
$\tilde\cE_{n+1}$ is local on $n$-cells it is possible to show by induction that each $u\in L^2(V_n)$ allows for a unique harmonic extension $\tilde u\in L^2(V_{n+1})$ and that 
\[
  \tilde\cE_{n+1}[\tilde u]=\frac{3}{5}\tilde\cE_n[u]
\]
for all $n\geq0$. This allows us to renormalise the bilinear form $\tilde\cE_n$ via 
\[
  \cE_{n}:=\left(\frac{5}{3}\right)^n\tilde\cE_n
\]
thus ensuring $\cE_{n+1}[\tilde u]=\cE_n[u]$ and therefore
$\cE_{n+1}[\hat u]\geq\cE_n[u]$ for any extension $\hat u$ of $u$.

Consider now any function $u:V^*\to\IR$ and denote its restriction to $V_n$ by $u_n$. Then, $\cE_n[u_n]$ is a non-decreasing sequence of nonnegative real numbers and therefore converges to
\[
  \cE[u]:=\lim_{n\to\infty} \cE_n[u_n]\in[0,\infty].
\]
Define $\scD(\cE)$ to be the set of all functions $u:V^*\to\IR$ for which this limit is finite. It can be shown that $u\in\scD(\cE)$ implies that $u$ is H\"older continuous and can therefore be uniquely extended to a continuous function on all of $\SG$. By abuse of notation we shall denote this extension by $u$ as well and set $\cE[u]:=\cE[u|_{V^*}]$ whenever the right-hand side is defined. Hence, by polarisation, we obtain a bilinear form $\cE(\cdot,\cdot)$ on $\scD(\cE)$. 

We further introduce the (probability) measure $\mu(\cdot):=\cH^s(\SG)^{-1}\cH^s(\cdot)$ on $\SG$, where $s:=\dimh(SG)=\frac{\log 3}{\log 2}$. It is possible to show that the bilinear form $(\cE,\scD(\cE))$ is a local regular Dirichlet form on $L^2(\SG,\mu)$ which in turn is attached to an operator as in equation \eqref{eq:DFtoOp}. This operator is what is known as standard Laplacian on $\SG$, and its spectral dimension is known to be $\frac{\log 3}{\log 5}$ (see e.g. \cite[section 3.5]{strichartz2006differential} or \cite{kigami1993weyl}). 

\subsubsection[Approximation by random walks]{Approximation by random walks\protect\footnote{The material in this section is an overview of the construction given in \cite[chapter II]{barlow1998diffusions}}}

It remains to discuss the walk dimension of the diffusion process generated by the standard Laplace operator on $\SG$. Once again, this construction uses the approximating graphs $G_n$. 

For $n\geq0$, let $Y^{(n)}_k, k\in\IN_0$ be a simple random walk on $G_n$, i.e. given a current position $Y^{(n)}_k=v\in V_n$, the process has equal probability to jump to each of the neighbours of $v$ in $V_n$. 

Take $v\in V_{n-1}$. Then, $v$ is contained in either one or two $(n-1)$-cells with conjunction points 
$\{x_1,...,x_i\}\ssq V_n\setminus\{v\}$, $i\in\{2,4\}$ depending on whether or not $v\in V_0$. Assume $Y^{(n)}_0=v$, $n\geq1$. The only way for $Y^{(n)}$ to leave the $(n-1)$-cells containing $v$ is via the points $\{x_1,...,x_i\}$. However, due to the symmetry, the probabilities of hitting $x_j, j=1,...,i$ first are equal. If we set 
\begin{align*}
  T^{n,m}_0&=\inf\left\{t\in\IN_0:Y^{(n)}_t\in V_m\right\}\\
  T^{n,m}_{k+1} &=\inf\left\{t\in\IN_0, t>T^{n+1,m}_k:
     Y^{(n)}_t\in V_m\setminus \left\{Y^{(n)}_{T^{n,m}_k}\right\}\right\},\ \ k\geq0,
\end{align*}
for $n>m$, then $\left(Y^{(n)}_{T^{n,m}_k}\right)_{k\in\IN_0}$ is a standard random walk on $G_m$ and therefore equal to 
$\left(Y^{(m)}_k\right)_{k\in\IN_0}$ in distribution. Using self-similarity and standard arguments for finite Markov-chains, it can be shown that $\PTE{T^{n,m}_{k+1}-T^{n,m}_k}=5^{n-m}$ and that $Y^{(n)}$ overcomes a distance (in graph metric) of $2^{n-m}$. The renormalised processes $X^{(n)}_k:=Y^{(n)}_{5^nk}$, $k\in\IN_0$ furthermore converge almost surely and uniformly on compact intervals to a continuous limit process $(X_t)_{t\geq0}$ with values in $\SG$. Additionally, the infinitesimal generator of $X$ coincides with the $\SG$-Laplacian introduced above. 
By construction, this process exhibits self-similarity, in the sense that in distribution, $(X_{5t})_{t\geq0}=(2X_t)_{t\geq0}$ for sufficiently small $t$, and it can indeed be verified that $\dimw(\SG,X)=\frac{\log 5}{\log 2}$

Putting everything together, we obtain for the Sierpi\'nski gasket:
\[
  \frac{\log 3}{\log 2}=\dimh(\SG)=\dims(\SG,\gD)\dimw(\SG,X)
  =\frac{\log 3}{\log 5}\frac{\log 5}{\log 2},
\]
so indeed, the Einstein relation holds on $\SG$.

\subsection{The real line as bounded metric space}\label{ssec:realline}

Bounded metric spaces form the most important class of spaces for which too naive of an adaption of \eqref{eq:defdimwrong} does not yield useful results. Indeed, consider the metric measure space 
$X=(\IR,d_{\arctan},\gl^1)$, where $d_{\arctan}(x,y)=|\arctan(x)-\arctan(y)|$ defines the metric. Since 
\[
  \tan:\left(\left(-\frac{\pi}{2},\frac{\pi}{2}\right),|\cdot|\right)\to(\IR,d_{\arctan}) 
\]
provides an isometry, we have $\dimh(X)=1$. On this space, we consider the negative of the usual weak Laplace operator, $-\gD_{\gl^1}$, defined by mapping a Sobolev function $u\in H_0^1(\IR,\gl^1)$ to the unique $g\in L^2(\IR,\gl^1)$ such that
\[
  \int_\IR g\gp\,\Di\gl^1=\int_\IR \partial_x u\partial_x \gp\,\Di\gl^1
\]
holds for all $\gp\in H_0^1(\IR,\gl^1)$. Notice how this does not differ from the negative weak Laplace operator on $(\IR,|\cdot|,\gl^1)$ since we did not change the measure and both metrics induce the same topology. Thus, we get from Weyl's classical result $\dims(X,-\gD_{\gl^1})=\frac{1}{2}$ and from the arguments developed in Section~6 that the Markov process associated to $\gD_{\gl^1}/2$ is $(B_t)_{t\geq0}$. 

It is now easy to see that \eqref{eq:defdimwrong} does not provide a useful notion of a walk dimension: Since $B_{\arctan}(x,R)=\IR$ for every radius $R\geq\pi$, the expression diverges to $\infty$. Even the more careful approach 
\[
  \lim_{R\nearrow\frac{\pi}{2}}\frac{\log \PTEp{\tau(R)}{0}}{\log R}
\]
runs into similar problems: Using the formula $\PTE{\tau([a,b])}=-ab$ for the exit time of a standard Brownian motion from the interval $[-a,b]\ni 0$, we get
\[
  \log \PTEp{\tau(B_{\arctan}(0,R))}{0}=2\log \tan R =\infty.
\]
However, the local walk dimension from Definition~\ref{def:dimw} works out quite elegantly: Setting $y=\arctan x\in\left(-\frac{\pi}{2},\frac{\pi}{2}\right)$, we obtain for some $\xi_1\in (y,y+r),\xi_2\in(y-r,y)$
\begin{multline*}
  \frac{\log \PTEp{\tau\big(B_{\arctan}(y,r)\big)}{y}}{\log r}
  =\frac{\log\left(\tan(y+r)-\tan y\right)}{\log r}+\frac{\log\left(\tan y-\tan (y-r)\right)}{\log r}\\
  =\frac{\log\frac{r}{\cos^2 \xi_1}}{\log r}+\frac{\log\frac{r}{\cos^2 \xi_2}}{\log r}
  =2-\frac{\log \cos^2 \xi_1}{\log r}-\frac{\log\cos^2 \xi_2}{\log r}
\end{multline*}
by using the mean value theorem. Taking the limit for $r\searrow0$ on both sides implies $\xi_1,\xi_2\to y$ and thus $\dimw(E,B_t,x)=2$.

\section{The Einstein Relation on Metric Measure Spaces}\label{sec:ERMMS}

This section is devoted to the investigation of the Einstein relation in the setting of an abstract mm-space. Let us recall that by an mm-space, we mean a complete separable metric space with a Radon measure. Whenever we want to be able to define the Einstein relation on an mm-space, we additionally assume that the space is locally compact, path-connected, and contains strictly more than one point. Furthermore, if the space is compact, we will always assume that the measure is a probability measure.

\subsection{The Einstein Relation under Lipschitz-isomorphisms}\label{ssec:Lipschitz}

\subsubsection{Lipschitz and mm-isomorphisms}

We will use this section to introduce two different categories $\mathsf{MM}_L$ and $\mathsf{MM}_{\leq1}$ whose objects are mm-spaces, but with different morphisms: 
\begin{itemize}
  \item In $\mathsf{MM}_L$, the set $\mathsf{MM}_L(X,Y)$ of morphisms from an object $X=(X,d_X,\mu_X)$ to another object $Y=(Y,d_Y,\mu_Y)$ is the set of all Lipschitz-continuous functions 
  \[ 
    \gp:\supp\mu_X\to\supp\mu_Y
  \]
  satisfying $\gp_*\mu_X=\mu_Y$, where $\gp_*\mu_X$ denotes the pushforward of $\mu_X$ by $\gp$.
  \item In $\mathsf{MM}_{\leq1}$, the set $\mathsf{MM}_{\leq1}(X,Y)$ of morphisms from an object $X=(X,d_X,\mu_X)$ to another object $Y=(Y,d_Y,\mu_Y)$ is the subset of $\mathsf{MM}_L(X,Y)$ consisting of all contraction maps (i.e. Lipschitz-continuous functions $f$ with $\Lip_f\leq1$, cf. \eqref{eq:Lip}).
\end{itemize}
In both of those categories, composition of morphisms is to be understood as the usual composition of maps. By definition, $\mathsf{MM}_{\leq1}$ is a subcategory of $\mathsf{MM}_L$. Considering the usual notion of isomorphism, both categories give rise to a meaningful concept of isomorphy for mm-spaces: 
\begin{definition}
  A Lipschitz-isomorphism between two mm-spaces $(X,d_X,\mu_X)$ and $(Y,d_Y,\mu_Y)$ is a map 
  $\gp:\supp\mu_X\to\supp\mu_Y$ with $\gp_*\mu_X=\mu_Y$ satisfying the bi-Lipschitz condition
  \[
    \frac{1}{C}d_X(x,y)\leq d_Y(\gp(x),\gp(y))\leq Cd_X(x,y)
  \]
  for all $x,y\in\supp\mu_X$ and a constant $C\in[1,\infty)$ not depending on $x,y$.
  
  Similarly, an mm-isomorphism is defined to be a Lipschitz-isomorphism with constant $C=1$. (This coincides with Definition~2.8 in \cite{shioya2016metric})
\end{definition}
As it turns out, Lipschitz-isomorphisms are precisely the isomorphisms in $\mathsf{MM}_L$, whereas mm-isomorphisms are the ones in $\mathsf{MM}_{\leq1}$.

Indeed, consider a Lipschitz-isomorphism $\gp:X\supseteq\supp\mu_X\to\supp\mu_Y\ssq Y$. By definition, this is an injective morphism from $\mathsf{MM}_L(X,Y)$. We need to show that $\gp$ is surjective to ensure the existence of a two-sided inverse in $\mathsf{MM}_L(Y,X)$. To this end, suppose there exists $y\in\supp\mu_Y\setminus \gp(\supp\mu_X)=:Z$. Since $\supp\mu_X$ is closed, so is its image under the homeomorphism $\gp$, and hence $Z\ssq \supp\mu_Y$ is open. As every open subset of $\supp\mu_Y$ is required to have positive measure, we obtain the contradiction
\[
  0<\mu_Y(Z)=\gp_*\mu_X(Z)=\mu_X\left(\gp^{-1}\left(\supp\mu_Y\setminus \gp(\supp\mu_X)\right)\right)=0.
\]
Hence, $\gp$ is indeed a bijection. Conversely, if $\gp$ is an isomorphism from $\mathsf{MM}_L(X,Y)$ then we get the lower bound from the Lipschitz-continuity of $\gp^{-1}\in\mathsf{MM}_L(Y,X)$, thus showing that $\gp$ is also a Lipschitz-isomorphism. Analogously, the corresponding statement for mm-isomorphisms can be derived.

We will write $(X,d_X,\mu_X)\simeq (Y,d_Y,\mu_Y)$ if $X$ and $Y$ are Lipschitz-isomorphic, whereas we will write 
$(X,d_X,\mu_X)\cong (Y,d_Y,\mu_Y)$ if they are mm-isomorphic. Trivially, $X\cong Y$ implies $X\simeq Y$. 

In what follows, we will always assume $\supp\mu_X=X$.
\begin{remark}
  We always have $(X,d_X,\mu_X)\cong(\supp\mu_X,d_X,\mu_X)$ by virtue of the identity map on $\supp\mu_X$. The restriction $\supp\mu_X=X$ becomes necessary for the Einstein relation since $\dimh(\supp\mu_X)$ might be strictly smaller than $\dimh(X)$, the term appearing in the Einstein relation \eqref{eq:ER}. We will later see (Proposition~\ref{prop:mmiso}) that the Einstein relation is invariant under Lipschitz-isomorphisms which provides some motivation to circumvent this restriction by considering the relation
  \[
    \dimh(\supp\mu_X)=c\dims(\supp\mu_X,A)\dimw(\supp\mu_X,M)
  \]
  instead of \eqref{eq:ER}.
\end{remark}

\subsubsection{Transport of structure}

Consider now two mm-spaces $(X,d_X,\mu_X)$ and $(Y,d_Y,\mu_Y)$ together with a map $\gp:X\to Y$, such that a suitable operator 
$A:L^2(X,\mu_X)\supseteq\scD(A)\to L^2(X,\mu_X)$ satisfies the Einstein relation with constant $c$ on $X$. How can we transport $A$ alongside $\gp$ to become an operator on $L^2(Y,\mu_Y)$, and which restrictions do we need to impose on $\gp$ to ensure that this transport of structure is compatible with the theory from Section~\ref{sec:Dims}?

Note first that any bimeasurable bijection $\gp:(X,d_X,\mu_X)\to(Y,d_Y,\mu_Y)$ induces by precomposition an operator
\begin{align*}
  \gp^*:L^2(Y,\mu_Y)&\to L^2(X,\mu_X)\\
  f&\mapsto f\circ\gp
\end{align*}
which is an isometric isomorphism because $\gp^{-1}:Y\to X$ induces its inverse $(\gp^*)^{-1} = (\gp^{-1})^*$ and because of
\begin{multline}\label{eq:isometry}
  \left\|\gp^*f\right\|_{L^2(X,\mu_X)}^2
  =\int_X |f(\gp(x))|^2\,\Di\mu_X
  =\int_Y |f|^2\,\Di\gp_*\mu_X\\
  =\int_Y |f|^2\,\Di\mu_Y
  =\|f\|_{L^2(Y,\mu_Y)}^2,
\end{multline}
by the change of variables formula for Lebesgue integrals. 

Denote by $\IL(H)$ the set of all partially defined linear maps (not necessarily bounded) on a Hilbert space $H$. Given an operator $A\in\IL(L^2(X,\mu_X))$, we can now construct an operator $\gp_{\IL} A\in\IL(L^2(Y,\mu_Y))$ by conjugating with $\gp^*$. More explicitly, we define the map
\[
  \gp_{\IL}:\IL(L^2(X,\mu_X))\to\IL(L^2(Y,\mu_Y))
\]
where $(\gp_{\IL} A)f:=((\gp^*)^{-1}\circ A\circ\gp^*)f$ and $\scD(\gp_{\IL} A)=(\gp^*)^{-1}(\scD(A))$. Note that $\gp_{\IL}$ is again a bijection with inverse given by $\gp_{\IL}^{-1}=(\gp^{-1})_{\IL}$ and that this bijection restricts to the spaces of bounded linear operators. 

It follows immediately that $\scD(\gp_{\IL} A)$ is dense if and only if $\scD(A)$ is, and that $\gp_{\IL} A$ is self-adjoint if and only if $A$ is. Indeed, consider arbitrary $f,g\in\scD(\gp_{\IL} A)$ with $f=(\gp^*)^{-1}(\bar f)$ and $g=(\gp^*)^{-1}(\bar g)$, where $\bar f, \bar g\in \scD(A)$. Then, applying \eqref{eq:isometry}, we have
\[
  \left<(\gp_{\IL} A)f,g\right>_{L^2(Y,\mu_Y)}
  =\left<(\gp^*)^{-1}A\gp^*(\gp^*)^{-1}\bar f,(\gp^*)^{-1}\bar g\right>_{L^2(Y,\mu_Y)}
  =\left<A\bar f,\bar g\right>_{L^2(X,\mu_X)}
\]
and we can perform the same calculations for $\left<f,(\gp_{\IL} A)g\right>_{L^2(Y,\mu_Y)}$, thus establishing the claimed equivalence. It is equally straightforward to check that the resolvent sets and the eigenvalues of $A$ and $\gp_{\IL} A$ coincide: Consider $\gl\in\rho(A)$, that is, 
$(\gl-A)^{-1}$ is a bounded linear operator on $L^2(X,\mu_X)$. To show $\gl\in\rho(\gp_{\IL} A)$, we consider
\begin{multline*}
  (\gl-\gp_{\IL} A)^{-1}
  =\left(\gl-(\gp^*)^{-1}A\gp^*\right)^{-1}
  =\left((\gp^*)^{-1}(\gl-A)\gp^*\right)^{-1}\\
  =(\gp^*)^{-1}(\gl-A)^{-1}\gp^*
  =\gp_{\IL}(\gl-A)^{-1}
\end{multline*}
which is a bounded linear operator on $L^2(Y,\mu_Y)$. If $\gl\in\IC$ is an eigenvalue of $A$ with eigenfunction $f\in L^2(X,\mu_X)$, then $(\gp^*)^{-1}f$ is an eigenfunction of $\gp_{\IL} A$ for the eigenvalue $\gl$ as well. This can easily be checked by calculating
\[
  (\gp_{\IL} A)\big((\gp^*)^{-1}f\big)=(\gp^*)^{-1}Af=\gl(\gp^*)^{-1}f.
\]

Moreover, $\gp_{\IL}$ respects operator semigroups: If $(T_t)_{t\geq0}$ is a strongly continuous contraction semigroup on $L^2(X,\mu)$ with generator $(-A,\scD(A))$ then $(\gp_{\IL} T_t)_{t\geq0}$ is a semigroup with the same properties on $L^2(Y,\mu_Y)$ and generator $(-\gp_{\IL} A,\scD(\gp_{\IL} A))$. Indeed, the semigroup property is trivial to check. For contractivity, note that for $L^2(Y,\mu_Y)\ni f=(\gp^*)^{-1}\bar f$ with $\bar f\in L^2(X,\mu_X)$,
\begin{multline*}
  \left\|(\gp_{\IL}T_t)f\right\|_{L^2(Y,\mu_Y)}
  =\left\|(\gp^*)^{-1}T_t\gp^*(\gp^*)^{-1}\bar f\right\|_{L^2(Y,\mu_Y)}
  =\left\|T_t\bar f\right\|_{L^2(X,\mu_X)}\\
  \leq\|\bar f\|_{L^2(X,\mu_X)}=\|f\|_{L^2(Y,\mu_Y)}.
\end{multline*}
For strong continuity, we calculate
\begin{multline*}
  \left\|\gp_{\IL} T_t \bar f-\gp_{\IL} T_0 \bar f\right\|_{L^2(Y,\mu_Y)}
  =\left\|(\gp^*)^{-1}(T_t\gp^*\bar f-\gp^*\bar f)\right\|_{L^2(Y,\mu_Y)}\\
  =\left\|T_t(\gp^*\bar f)-(\gp^*\bar f)\right\|_{L^2(X,\mu_X)} \to 0
\end{multline*}
for $t\searrow 0$ and arbitrary $\bar f\in L^2(X,\mu)$, and verifying the generator works analogously. 

Note however that a bi-measurable bijection $\gp$ does not respect enough structure to ensure that 
$\gp_\IL \sqrt{A}$ generates a regular Dirichlet form if and only if $\sqrt{A}$ does -- recall that this means the density of $\scD(\gp_\IL\sqrt{A})\cap C_c(Y)$ in both $\scD(\gp_\IL\sqrt{A})$ and $C_c(Y)$. To this end, suppose now that $\gp:X\to Y$ is a homeomorphism between $X$ and $Y$ (since both spaces are equipped with their Borel $\gs$-algebras, such $\gp$ is automatically bi-measurable and bijective). Similar to the case of $L^2$-spaces, this induces an isometric isomorphism $\gp^*:C_0(Y)\to C_0(X), \gp^*(f)=f\circ\gp$ between algebras of continuous functions vanishing at infinity, equipped with sup-norm $\|\cdot\|_{C_0}$. This isomorphism restricts to the subalgebras of compactly supported continuous functions $C_c(X)$ resp. $C_c(Y)$. 

\begin{lemma}
  With the notation just introduced, if the Dirichlet form on $L^2(X,\mu_X)$ defined by 
  \[
    \cE(f,g):=\left<\sqrt{A}f,\sqrt{A}g\right>_{L^2(X,\mu_X)}\ \text{ for } f,g\in\scD(\sqrt{A})
  \]
 is regular then so is the Dirichlet form on $L^2(Y,\mu_Y)$ defined by 
 \[
   (\gp_\IL)^*\cE(\bar f,\bar g):=\left<(\gp_\IL\sqrt{A})\bar f,(\gp_\IL\sqrt{A})\bar g\right>_{L^2(Y,\mu_Y)}\ \text{ for } \bar f,\bar g\in\gp_*^{-1}(\scD(\sqrt{A})).
 \]
\end{lemma}
\begin{proof}
  We need to show that the intersection of $\scD(\gp_\IL\sqrt{A})=(\gp^*)^{-1}(\scD(\sqrt{A}))$ and $C_c(Y)$ is dense both in $C_c(Y)$ with respect to $\|\cdot\|_{C_0(Y)}$ and in $\scD(\gp_\IL\sqrt{A})$ with respect to $((\gp_\IL)^*\cE)_1$ as introduced in Definition~\ref{defin:DF}.
  
  For the first part, take $C_c(Y)\ni f=(\gp^*)^{-1} g$ for $g\in C_c(X)$. Then, there exists a sequence $(g_n)_{n\in\IN}\ssq \scD(\sqrt{A})\cap C_c(X)$ with $\|g_n-g\|_{C_0(X)}\to 0$ as $n\to\infty$. Since $(\gp^*)^{-1}$ is isometric, we conclude that $f_n:=(\gp^*)^{-1}g_n\in \scD(\gp_\IL\sqrt{A})\cap C_c(Y)$ converges to $f$ in 
  $\|\cdot\|_{C_0(Y)}$.
  
  For the second part, we take $\scD(\gp_\IL\sqrt{A})\ni f=(\gp^*)^{-1} g$ for $g\in\scD(\sqrt{A})$. By regularity of $\cE$, there exists again a sequence $(g_n)_{n\in\IN}\ssq\scD(\sqrt{A})\cap C_c(X)$ with $\cE_1[g_n-g]\to 0$ as $n\to\infty$. Setting $f_n:=(\gp^*)^{-1} g_n\in\scD(\gp_\IL\sqrt{A})\cap C_c(Y)$ and recalling the notation from Definition~\ref{defin:DF}, we obtain
  \begin{align*}
    ((\gp_\IL)^*&\cE)_1[f_n-f]
     =\left(\left<\gp_\IL\sqrt{A}\ \cdot\ ,\gp_\IL\sqrt{A}\ \cdot\right>_{L^2(Y,\mu_Y)} + \left<\cdot,\cdot\right>\right)[f_n-f]\\
    &=\left\|\left(\gp_\IL\sqrt{A}\right)(f_n-f)\right\|_{L^2(Y,\mu_Y)}^2
        +\|f_n-f\|_{L^2(Y,\mu_Y)}^2\\
    &=\left\|(\gp^*)^{-1}\sqrt{A}\gp^*(\gp^*)^{-1}(g_n-g)\right\|_{L^2(Y,\mu_Y)}^2
        +\|(\gp^*)^{-1}(g_n-g)\|_{L^2(Y,\mu_Y)}^2\\
    &=\left\|\sqrt{A}(g_n-g)\right\|_{L^2(X,\mu_X)}^2+\|g_n-g\|_{L^2(X,\mu_X)}^2\\
    &=\cE_1[g_n-g]\to 0,
  \end{align*}
  which concludes the proof.
\end{proof}

Putting everything together, we observe that $A$ satisfies the assumptions in \ref{cond:A} if and only if $\gp_\IL A$ does whenever $\gp:X\to Y$ is a homeomorphism, and then $\dims(X,A)=\dims(Y,\gp_\IL A)$. The spectral dimension is therefore stable under a very large class of transformations. As it turns out, this will not be the case for Hausdorff and walk dimension. 

\begin{proposition}\label{prop:mmiso}
  Let $(X,d_X,\mu_X)$ and $(Y,d_Y,\mu_Y)$ be complete separable locally compact path-connected metric measure spaces with $\supp\mu_X=X$ and $\supp\mu_Y=Y$ that are Lipschitz-isomorphic by virtue of the map $\gp:X\to Y$. Suppose the Einstein relation with constant $c$ holds on $X$ with respect to an operator $(A,\scD(A))$ satisfying assumptions \ref{cond:A}. Then, the Einstein relation also holds on $Y$ with the same constant $c$ and with respect to $\gp_{\IL} A$.
\end{proposition}
\begin{proof}
  As the Hausdorff dimension is invariant under bi-Lipschitz maps we obtain $\dimh(X)=\dimh(Y)$, and as observed above, $\dims(X,A)=\dims(Y,\gp_{\IL} A)$. So, it remains to show $\dimw(X,M)=\dimw(X,M^{(\gp)})$ where $M$ is a Hunt process associated to $A$ and $M^{(\gp)}$ is one associated to $\gp_{\IL} A$. 
  
  We consider the process $N_t:=\gp(M_t)$. This process is a Hunt process with values in $Y$, and possesses the semigroup
  \begin{multline*}
    T_t^{(N)}f
    =\PTEp{f(N_t)}{\cdot}
    =\PTEp{(f\circ\gp)(M_t)}{\gp^{-1}(\cdot)}\\
    =T_t[\gp^*f](\gp^{-1}(\cdot))
    =(\gp^*)^{-1}T_t\gp^* f
    =(\gp_{\IL} T_t)f
  \end{multline*}
  where we used the notation from the discussion above. 
  
  Thus, due to Theorem~\ref{thm:fukushima}, the processes $N$ and $M^{(\gp)}$ coincide up to their behaviour on a polar set. It is therefore enough to determine the walk dimension for $N_t$. By the bi-Lipschitz continuity of $\gp$, we obtain 
  \[
    \gp\left(B_X\left(x,C^{-1}r\right)\right)
    \ssq B_Y\big(\gp(x),r\big) 
    \ssq \gp\big(B_X(x,Cr)\big)
  \]
  where $C>0$ is the two-sided Lipschitz constant of $\gp$. Hence, 
  \[
   \tau_M(C^{-1}r)\leq \tau_N(r)\leq \tau_M(Cr).
  \]
  From this, we get for all sufficiently small $r>0$
  \[
    \frac{\log Cr}{\log r}
     \cdot\frac{\log \PTEp{\tau_M(Cr)}{x}}{\log Cr}
    \leq \frac{\log \PTEp{\tau_N(r)}{\gp(x)}}{\log r}
    \leq \frac{\log C^{-1}r}{\log r}
     \cdot\frac{\log \PTEp{\tau_M(C^{-1}r)}{x}}{\log C^{-1}r}.
  \]
  Taking the limit for $r\searrow0$ and applying a standard squeezing argument, we obtain $\dimw(X,M)=\dimw(Y,N)$. 
\end{proof}

\begin{remark}\label{rem:ToS}
  Note that we required the bi-Lipschitz property for determining 
  $\dimh$ and $\dimw$, whereas we only needed $\gp$ to be a homeomorphism in order to show that $M^{(\gp)}$ and $N$ share the same semigroup. This allows us in the following sections -- given a homeomorphism $\gp:(X,d_X,\mu_X)\to(Y,d_Y)$ -- to transport the complete structure needed for the Einstein relation by 
  \begin{itemize}
    \item Endowing $(Y,\mu_Y)$ with the push-forward measure $\gp_*\mu_X$.
    \item Mapping the generator $(A,\scD(A))$ to 
    $(\gp_\IL A,(\gp^*)^{-1}\scD(A))$, thus also mapping the generated semigroup $T_t$ to $\gp_\IL T_t$.
    \item Sending the Hunt process $M$ to $\gp(M)$. 
  \end{itemize}
  What we did so far ensures that all these constructions are compatible with each other. 
\end{remark}

From Proposition~\ref{prop:mmiso} we immediately obtain the following two corollaries:
\begin{corollary}
  If $(X,d_X,\mu_X)\cong(Y,d_Y,\mu_Y)$ and the Einstein relation with constant $c$ holds on $X$ w.r.t. $(A,\scD(A))$ is an operator on $L^2(X,\mu_X)$, then it also holds on $Y$ with the same constant w.r.t. $\gp_{\IL}A$. 
\end{corollary}
\begin{corollary}
  If $X\ssq\IR^n$ and $d_1$ and $d_2$ are metrics which are induced by norms, then $\id_X:(X,d_1,\mu)\to(X,d_2,\mu)$ will preserve the constant in the Einstein relation.
\end{corollary}
The second corollary follows from the well-known fact that all norms on a finite-dimensional Banach space are equivalent.

\subsection{H\"older regularity and graphs of functions}\label{ssec:Holder}

A natural question arising at this point is whether the invariance of the Einstein relation of Proposition~\ref{prop:mmiso} can be extended to a larger class of transformations. In particular, what happens if $\gp$ is only a H\"older continuous map instead of a bi-Lipschitz one?

As we saw in the previous subsection, such $\gp$ does not impede the spectral dimension, but it is well-known that $\ga$-H\"older continuous transformations are not compatible with the Hausdorff dimension, besides the general estimate 
\[
 \dimh(\gp(X))\leq\ga^{-1}\dimh(X)
\]
mentioned in Subsection~\ref{ssec:Hausdorff}. We will see that a similar picture occurs for the walk dimension.

\begin{definition}
  Let $\ga\in(0,1]$. We say that a map $\gp:(X,d_X)\to(Y,d_Y)$ between two metric spaces is locally $\ga$-H\"older continuous at $x\in X$ if there exists an open neighbourhood $U\ssq X$ of $x$ and a constant $C>0$ such that
  \[
    d_Y(\gp(x),\gp(y))\leq Cd_X(x,y)^\ga
  \]
  for all $y\in U$. If this holds for all $x\in X$ we call $\gp$ locally $\ga$-H\"older continuous on $X$. 
\end{definition}
Note that if $\gp$ is $\ga$-H\"older continuous then it is also 
$\gb$-H\"older continuous for any $\gb<\ga$ and that for $\ga=1$, we get back the definition of Lipschitz continuity. This allows us to define H\"older regularity as precisely the parameter $\ga$ at which the phase transition between being H\"older continuous and not being H\"older continuous occurs.
\begin{definition}
  In extension of the previous definition, we say that $\gp$ is 
  locally $\ga$-H\"older regular at $x\in X$ if $\ga$ is the supremum of all $\gb>0$ for which $\gp$ is locally $\gb$-H\"older continuous at $x$. Equivalently, such $\ga$ is the infimum of $1$ and all 
  $\gc\leq1$ for which $\gp$ is not locally $\gc$-H\"older continuous at $x$. 
\end{definition}

\subsubsection{A closer look at the walk dimension}

\begin{lemma}\label{lem:dimuw1}
  Let $\left(M_t^x\right)_{t\geq0}$ be a right-continuous stochastic process on $(X,d_X)$ starting in $x\in X$ and let $\gp:(X,d_X)\to(Y,d_Y)$ be a map which is locally $\ga$-H\"older regular at $x$ for some $\ga\in(0,1)$. Suppose further that the local walk dimension of $M$ at $x$ exists. Then the upper local walk dimension, defined by
  \[
    \dimuw(X,M;x):=\limsup_{r\searrow0}
       \frac{\log\PTEp{\tau_M(r)}{x}}{\log r},
  \]
  satisfies 
  \begin{equation}\label{eq:estdimuw}
    \dimuw\big(Y,\gp(M);\gp(x)\big)
      \leq\frac{1}{\ga}\dimw(X,M;x).
  \end{equation}
\end{lemma}
\begin{proof}
  Let $0<\gb<\ga$. Then, $\gp$ is locally $\gb$-H\"older continuous at $x$ and therefore, there exists a constant $C>0$ such that 
  \[
    \gp\left(B_X(x,r)\right)\ssq B_Y\left(\gp(x),Cr^\gb\right) 
  \]
  for all sufficiently small $r>0$. Thus, if $\gp(M)$ exits 
  $B_Y\left(\gp(x),Cr^\gb\right)$, it already left $\gp(B_X(x,r))$.
  Since $\gp(M_t)\notin\gp(B_X(x,r))$ implies $M_t\notin B_X(x,r)$ by comparing the preimages, we obtain the inequality
  \[
    \tau_M(r)
    \leq \tau_{\gp(M)}\big(\gp(B_X(x,r))\big)
    \leq \tau_{\gp(M)}(Cr^\gb).
  \]
  This implies for all $r$ small enough
  \[
    \frac{\log\PTEp{\tau_{\gp(M)}(Cr^\gb)}{\gp(x)}}{\log Cr^\gb}
    \leq \frac{\log\PTEp{\tau_M(r)}{x}}{\log r}
      \cdot \frac{\log r}{\log Cr^\gb},
  \]
  where the right-hand side converges to $\gb^{-1}\dimw(X,M;x)$ as $r\searrow0$, thus showing
  \[
    \dimuw\big(Y,\gp(M);\gp(x)\big)\leq \frac{1}{\gb}\dimw(X,M;x).
  \]
  As all estimates are valid for every $\gb<\ga$ and since
  $\dimuw\big(Y,\gp(M);\gp(x)\big)$ does not depend on $\gb$, we can take the supremum over all $\gb<\ga$ to obtain
  \[
    \dimuw\big(Y,\gp(M);\gp(x)\big)\leq \frac{1}{\ga}\dimw(X,M;x)
  \]
  which concludes the proof.
\end{proof}
In general, equality in \eqref{eq:estdimuw} does not hold. Fix $0<\ga<1$. We consider the measure space $(\IR^2,\gl^2)$ endowed with two different metrics -- first with the metric 
\[
  d_1^{(\ga)}(x,y)=|x_1-y_1|^\ga+|x_2-y_2|
  \ \ \text{ for }\ x=(x_1,x_2),y=(y_1,y_2)\in\IR^2
\]
and second with the metric $d_1$ induced by the 1-norm. In other words, we set $X=\left(\IR^2,d_1^{(\ga)},\gl^2\right)$ and 
$Y=\left(\IR^2,d_1,\gl^2\right)$. By definition, 
$\id:X\to Y$ provides a homeomorphism that is everywhere locally 
$\ga$-H\"older continuous. Let $(W_t)_{t\geq0}$ be a 1-dimensional standard Wiener process, and regard $(0,W_t)$ as a process in $X$ which has $\dimw(X,(0,W_t),0)=2$. Therefore,
\[
  \dimw(Y,(0,W_t),0)=2<\frac{2}{\ga}.
\]
Despite this counterexample, we get equality in \eqref{eq:estdimuw} in the following setting:
\begin{lemma}\label{lem:dimuw2}
  Let $(X,d_X)$ be a path-connected metric space consisting of more than one single point and let $M=(M_t^x)_{t\geq0}$ be an $X$-valued continuous stochastic process starting in $x\in X$. Let 
  $\gp:(X,d_X)\to(Y,d_Y)$ be locally $\ga$-H\"older regular at $x$. Suppose further that $\dimw(X,M;x)$ exists. Then
  \[
    \dimuw\big(Y,\gp(M);\gp(x)\big)=\frac{1}{\ga}\dimw(X,M;x)
  \]
  holds, provided there exists a constant $C>1$ and a sequence 
  $(r_n)_{n\in\IN}$ with $r_n\searrow 0$ such that for all $n\in\IN$, there exists a set $\gC_n\ssq B_X(x,Cr_n)\setminus B_X(x,r_n)$ subject to the following two conditions:
  \begin{enumerate}
    \item For all $y\in\gC_n$, $\gp$ violates an $\ga$-H\"older estimate: 
    $d_Y(\gp(x),\gp(y))>d_X(x,y)^{\ga+\gep}$.
    \item The complement of $\gC_n$, $X\setminus\gC_n$, splits into at least two non-empty path-connected components. 
  \end{enumerate}
\end{lemma}
\begin{proof}
  Due to the previous lemma, it only remains to show ``$\geq$''. By 
  $X_n$, we denote the connected component of $X\setminus\gC_n$ which contains $x$. Since we have the inclusions $B_X(x,r_n)\ssq X_n \ssq B_X(x,Cr_n)$ by definition of $\gC_n$, assumption $1.$ implies that
  \[
    B_Y\left(\gp(x),r_n^{\ga}\right)
    \ssq Y_n:=\gp(X_n)
    \ssq \gp\left(B_X(x,Cr_n)\right).
  \]
  Thus, 
  \[
    \tau_{\gp(M)}\left(r_n^{\ga}\right)
    \leq \tau_{\gp(M)}(Y_n)
    \leq \tau_{\gp(M)}\left(\gp\left(B_X(x,Cr_n)\right)\right)
    =\tau_M(Cr_n),
  \]
  which in turn yields to
  \[
    \frac{\log\PTEp{\tau_M(Cr_n)}{x}}{\log Cr_n}
      \cdot\frac{\log Cr_n}{\log r_n^{\ga}}
    \leq \frac{\log\PTEp{\tau_{\gp(M)}\left(r_n^{\ga}\right)}{\gp(x)}}
      {\log r_n^{\ga}}
  \]
  and consequentially for $n\to\infty$ to
  \begin{equation}\label{eq:w}
    \frac{1}{\ga}\dimw(X,M;x)\leq \liminf_{n\to\infty}
    \frac{\log\PTEp{\tau_{\gp(M^x)}\left(r_n^{\ga}\right)}{\gp(x)}}
      {\log r_n^{\ga}}=: w.
  \end{equation}
  By virtue of Lemma~\ref{lem:dimuw1} we also have 
  \[
    w\leq \dimuw(Y,\gp(M);\gp(x))\leq \frac{1}{\ga}\dimw(X,M;x)
  \]
  which shows the assertion when combined with \eqref{eq:w}.
\end{proof}
\begin{remark} 
  Suppose $N=\left(N_t^y\right)_{t\geq0}$ is a stochastic process with values in $(Y,d_Y)$ that is almost surely (locally) $\ga$-H\"older regular and satisfies the assumptions of Lemma~\ref{lem:dimuw2} with probability 1. Then we can always use Lemma~\ref{lem:dimuw2} to obtain
  $\dimuw(Y,N;\cdot\ )=\ga^{-1}$ by setting 
  $(X,d_X)=(\IR_{\geq0},|\cdot|)$, choosing $M$ deterministically as $M_t=t$ and regard $N$ as the (random) map $\gp=N:X\to Y$.
\end{remark}
In the special case where $X$ is an open domain in the 1-dimensional euclidean space and $M$ is a Brownian motion in $X$ we can disregard condition $2.$ in Lemma~\ref{lem:dimuw2} since the exit time for the Brownian motion does only depend on the distance from the starting point. 

\subsubsection{Graphs of continuous functions}

Given a continuous map $f:(X,d_X)\to (Y,d_Y)$ between two metric spaces, its graph
\[
  \gr(f):=\left\{(x,f(x))\in X\times Y:x\in X\right\}
\]
can be equipped with the restriction of the maximum metric on 
$X\times Y$,
\[
  d_\infty\big((x,y),(x',y')\big):=d_X(x,x')\vee d_Y(y,y'),\ \ 
  x,x'\in X, y,y'\in Y,
\]
to $\gr(f)\ssq X\times Y$. This makes $(\gr(f),d_\infty)$ a metric space that comes with a natural map $\gp:X\to\gr(f)$ sending $x\in X$ to 
$\gp(x):=(x,f(x))$. Since $f$ is continuous, it is easy to check that 
$\gp$ provides a homeomorphism between $X$ and $\gr(f)$ with the inverse given by the projection onto the first coordinate, $\pi$.
We point out that while $\pi$ is always Lipschitz-continuous, $\gp$ is 
(locally) $\ga$-H\"older continuous if and only if $f$ is. Indeed, we have 
\begin{multline*}
  d_Y(f(x),f(x'))<Cd_X(x,x')^\ga \\
  \Longleftrightarrow
  d_X(x,x')\vee d_Y(f(x),f(x'))<(1\vee C)d_X(x,x')^\ga
\end{multline*}
whenever $d_X(x,x')<1$. 

This setting is therefore a natural application to the arguments of the previous section. Unfortunately, not much is known about the Hausdorff dimension of these objects.

As deterministic $\ga$-H\"older regular functions are rather complicated objects to deal with, we will instead consider random functions. More precisely, we will look at 1-dimensional continuous $\ga$-self-similar process $(X_t)_{t\in\IR}$ with stationary increments over a suitable probability space $(\gO,\scA,\Prob)$. Here, $\ga$-self-similar for 
$0<\ga<1$ means that the processes $(X_t)_{t\in\IR}$ and 
$\left(\xi^{-\ga}X_{\xi t}\right)_{t\in\IR}$ have the same distribution for any $\xi>0$. By a theorem of Taqqu, see \cite[Theorem 1.3.1]{embrechts2002selfsimilar}, such a process is automatically a fractional Brownian motion, up to a constant factor.  

Fix $H\in(0,1)$ and recall that a fractional Brownian motion $B^H=\left(B^H_t\right)_{t\in\IR}$ with Hurst index $H$ is the centered Gaussian process with $B^H_0=0$ and covariance function
\[
  \PTE{B^H_s,B^H_t}=\frac{1}{2}\left(t^{2H}+s^{2H}-|t-s|^{2H}\right)\quad \text{ for }s,t\geq0.
\]
It is easy to check that this defines a $H$-self-similar process. By the Theorems~4.1.1 and 4.1.3 in \cite{embrechts2002selfsimilar}, there exists a version of $B^H$ which is almost surely everywhere locally $H$-H\"older regular (note though that $B_H$ is only $\ga$-H\"older continuous for $\ga<H$, and not for $\ga=H$). Note also that we consider a 2-sided fractional Brownian motion, to avoid boundary issues at $t=0$. 

As pointed out in Remark~\ref{rem:ToS}, we can now transfer the analytic structure of Section~6 on the real line $\IR$ via $\gp$ to 
$\gr(B^H)$. More explicitly, we have the measure $\gp_*\gl^1$ on $\gr(B^H)$ and an operator $\gp_L\gD_{\gl^1}$ acting on $L^2(\gr(B^H),\gp_*\gl^1)$ that generates the Hunt process $(\gp(W_s^t))_{s\geq0}$, where $(W_s^t)_{s\geq0}$ is a Wiener process independent from $B^H$ with start in $t\in\IR$.

We note the following:
\begin{itemize}
  \item From \cite{adler1977hausdorff}, we get $\dimh(\gr(B^H))=2-H$ with probability 1.  
  \item As discussed in the last section we have
  \[
    \frac{1}{2}=\dims\left(\IR,-\frac{1}{2}\gD_{\gl^1}\right)
    =\dims\left(\gr(B^H),\Phi \left(-\frac{1}{2}\gD_{\gl^1}\right)\right),
  \]
  where we once again appealed to Weyl's classical results. 
  \item The walk dimension is given by $\dimw(\gr(B^H),\gp(W))=\frac{2}{H}$ with probability 1, that is, we have the next theorem:
\end{itemize}

\begin{theorem}\label{thm:BHdimw}
 Let $B^H=\left(B_t^H\right)_{t\in\IR}$ be a 2-sided fractional Brownian motion on with Hurst index $H\in (0,1)$, and denote by $\gr(B^H)\ssq\IR^2$ its graph. Furthermore, let $(W_t)_{t\geq0}$ be a Wiener process independent from $B^H$, and denote by $\gp:\IR\to\gr(B^H)$ the map $x\mapsto (x,B^H_x)$. Then, with probability 1, the walk dimension $\dimw(\gr(B^H),\gp(W))$ exists and is equal to $\frac{2}{H}$.
\end{theorem}

Thus, the Einstein relation, despite holding with constant 1 on 
$\IR$ with $-\frac{1}{2}\gD_{\gl^1}$, changes its constant under application of $\gp$ to
\[
  c(H)=H(2-H),\ \ H\in(0,1).
\]
This is remarkable, as we generally only have the upper bound $2-\ga$ for both $\dimh$ and $\dimw$ under $\ga$-H\"older regular transformations $\IR\to\IR$ -- cf. \cite[chapter 16]{falconer2007fractal} for the upper bound on the Hausdorff dimension. The case exhibited here provides an example where both dimensions are changed differently.

For the proof of Theorem~\ref{thm:BHdimw}, observe first that according to Lemma~\ref{lem:dimuw1} we obtain $2/H$ as upper bound to the upper local walk dimension at any point of $\gr(B^H)$. In what follows, it is therefore enough to show that $2/H$ is also a lower bound to the walk dimension.

Consider a point $x=x(\go)=(T,B^H_T(\go))\in\gr(B^H_\cdot(\go))$ and chose $r>0$. Let $B_\infty(x,r)$ denote the open ball of radius $r$ around $x$ with respect to $d_\infty$. Introduce further the random times $\gT^+_r\big(B^H_\cdot(\go),T\big)$ and $\gT^-_r\big(B^H_\cdot(\go),T\big)$, denoting the time where the process $B^H$ first exits resp. last enters $B_\infty(x,r)$ -- in other words,
\begin{align}\label{eq:Theta}
 \gT^+_r\big(B^H_\cdot(\go),T\big)
 &:=\inf\left\{t>T:(t,B^H_t)\notin B_\infty(x,r)\right\}\notag\\
    \gT^-_r\big(B^H_\cdot(\go),T\big)
 &:=\sup\left\{t<T:(t,B^H_t)\notin B_\infty(x,r)\right\}.
\end{align}
By the standard result for the expectation of two-sided exit times for the Wiener process, we now obtain 
\[
 \PTEp{\tau_{\gp(W)}(B_\infty(x,r))}{x}
 =-\Big(\gT^+_r\big(B^H_\cdot(\go),T\big)-T\Big)
   \Big(\gT^-_r\big(B^H_\cdot(\go),T\big)-T\Big)
\]    
and consequentially
\begin{equation}\label{eq:sumoflogs}
 \frac{\log \PTEp{\tau_{\gp(W)}(B_\infty(x,r))}{x}}{\log r}
 =\frac{\log \left(\gT^+_r\big(B^H_\cdot(\go),T\big)-T\right)}{\log r}
  +\frac{\log \left(T-\gT^-_r\big(B^H_\cdot(\go),T\big)\right)}{\log r}.
\end{equation}
We will show that the limit inferior of each summand on the right-hand side is bounded from below by $1/H$ as $r\searrow 0$. To this end, we further introduce for $h^\pm\geq0$ the random variables $\gvt_{r,h^-,h^+}^+(T)$ and $\gvt_{r,h^-,h^+}^-(T)$ by 
\begin{align}\label{eq:gvtdef}
 \gvt_{r,h^-,h^+}^\pm(T)
 :&=\inf\left\{s\geq0:B^H_{T\pm s}\notin(B^H_T-r-h^-,B^H_T+r+h^+)\right\}\notag\\
  &=\inf\left\{s\geq0:B^H_{T\pm s}-B^H_T\notin(-r-h^-,r+h^+)\right\}.
\end{align}
These random variables relate to \eqref{eq:sumoflogs} via
\begin{equation*}
 \gvt_{r,0,0}^\pm(T) \wedge r=\left|\gT_r^\pm(B^H,T)-T\right|.
\end{equation*}

First, we show the following technical lemma:
\begin{lemma}\label{lem:BMdimuw}
 Assume that $h^\pm=h^\pm(r)$ are nonnegative functions defined on some interval $(0,r_0)$ such that $h^\pm=O(r)$ as $r\searrow0$. For fixed $T_0\in\IR$, we then have 
 \begin{equation*}
  \lim_{r\searrow0}\frac{\log \gvt_{r,h^-,h^+}^\pm(T_0)}{\log r} =\frac{1}{H}
 \end{equation*}
 $\Prob$-almost surely.
\end{lemma}
\begin{proof}
 From the definition of $\gvt^\pm$ and the self-similarity and the stationary increments of $B^H$, we obtain
 \begin{align*}
  \gvt_{\xi r,h^-,h^+}^\pm(T_0)
  &=\inf\left\{s\geq0:B^H_{T_0\pm s}-B^H_{T_0} \notin(-\xi r-h^-,\xi  r+h^+)\right\}\\
  &\stackrel{\cD}{=}\inf\left\{s\geq0:B^H_{\pm s}\notin
    (-\xi r-h^-,\xi r+h^+)\right\}\\
  &=\inf\left\{s\geq0:\xi^{-1}B^H_{\pm s}\notin(-r-h^-/\xi,r+h^+/\xi)\right\}\\
  &\stackrel{\cD}{=}\inf\left\{s\geq0:B^H_{\pm\xi^{-1/H}s}
    \notin (-r-h^-/\xi,r+h^+/\xi)\right\}\\
  &=\xi^{1/H}\inf\left\{s\geq0:B^H_{\pm s}\notin(-r-h^-/\xi,r+h^+/\xi)\right\}\\
  &\stackrel{\cD}{=}\xi^{1/H}\gvt_{r,h^-/\xi,h^+/\xi}^{\pm}(T_0)
  \end{align*}
  for any $\xi>0$. Hence we have 
  \begin{equation}\label{eq:gvtss}
   \gvt_{r,h^-,h^+}^\pm(T_0)\stackrel{\cD}{=}r^{1/H}\gvt_{1,h^-/r,h^+/r}^\pm(T_0).
  \end{equation}
  Observe further that $\gvt_{r,h^-,h^+}^\pm(T_0)$ is monotonous in both $h^-$ and $h^+$. Thus by assumption, there exists a positive constant $C$ not depending on $r$ such that 
  \[
   \gvt_{1,0,0}^\pm(T_0)\leq \gvt_{1,h^-/r,h^+/r}^\pm(T_0) \leq \gvt_{1,C,C}^\pm(T_0)
  \]
  holds for all $r>0$ sufficiently small. Together with \eqref{eq:gvtss}, this yields
  \[
   \frac{1}{H}+\frac{\log \gvt_{1,C,C}^\pm(T_0)}{\log r} \leq \frac{\log \gvt_{r,h^-,h^+}^\pm(T_0)}{\log r} \leq \frac{1}{H}+\frac{\log \gvt_{1,0,0}^\pm(T_0)}{\log r}.
  \]
  Here, the left- and right-hand sides converge to $1/H$ as $r\searrow0$ since $\Prob$-almost surely, $0<\gvt_{1,0,0}^\pm(T_0),\gvt_{1,C,C}^\pm(T_0)<\infty$ holds. 
\end{proof}

\begin{figure}[ht]
\centering
\includegraphics[width=0.95\textwidth]{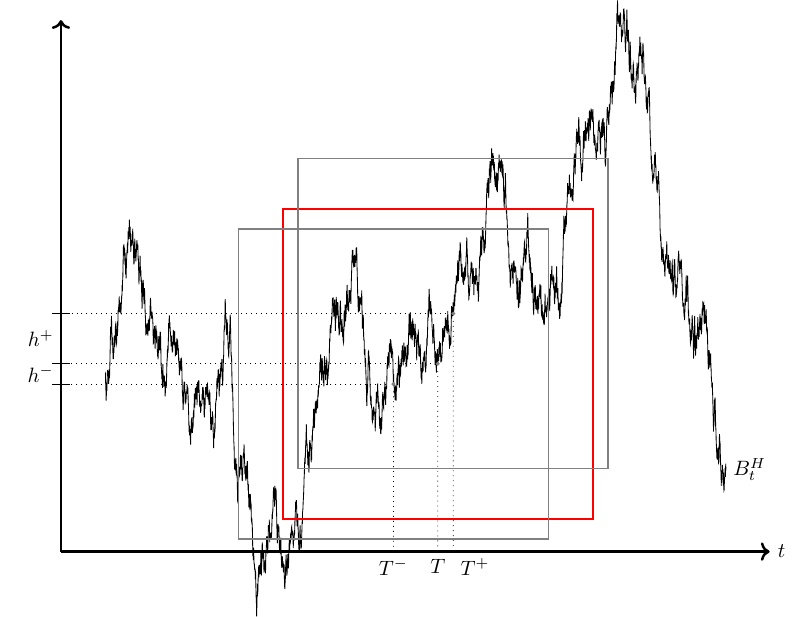}
\caption{A fractional Brownian motion, with the balls $B_\infty(\cdot,r)$ indicated for the points over $T$ (in red) and $T^\pm$ (in grey) on the graph.}\label{fig:bm}
\end{figure}

As an immediate consequence, we note that the claim of Lemma~\ref{lem:BMdimuw} holds $\Prob$-almost surely for all $T_0\in\IQ$ simultaneously (where we can even allow $h^\pm$ to depend on $T_0$).

Consider now again an arbitrary $T\in\IR$, and assume that $r<1$ is sufficiently small for $\left|B^H_{T'}(\go)-B^H_T(\go)\right|\leq \left|T'-T\right|^{H/2}$ to hold whenever $|T'-T|\leq r$. For $T$, choose rational numbers $T^\pm\in(T-r^{2/H},T+r^{2/H})$ with $T^-<T<T^+$. Moreover, set $h^\pm=\left|B^H_T(\go)-B^H_{T^\pm}(\go)\right|$, and observe that this implies $h^\pm\leq|T^\pm-T|^{H/2}\leq r$. In particular, as functions in $r$, we have $h^\pm(r)=O(r)$.

We will now bound $\big|\gT_r^\pm(B^H,T)-T\big|$ in terms of $\gvt^\pm$ at the nearby points $T^\pm$, see Figure~\ref{fig:bm} for a sketch. Indeed, $B_\infty(x,r)$ is a subset of both horizontal ``stripes'' $\IR\times\big(B^H_{T^\pm}-r-h^\pm,B^H_{T^\pm}+r+h^\pm\big)$, and therefore
\begin{align*}
 \left|\gT_r^\pm(B^H,T)-T\right|
 &\leq \max\{|T^--T|,|T^+-T|\} + \max\Big\{\gvt_{r,h^-,h^-}^\pm(T^-),\gvt_{r,h^+,h^+}^\pm(T^+)\Big\}\\
 &\leq r^{4/H}+r^{1/H+o(1)}=r^{1/H+o(1)},
\end{align*}
where we applied Lemma~\ref{lem:BMdimuw} to $\gvt_{r,h^-,h^-}^\pm(T^-)$ and $\gvt_{r,h^+,h^+}^\pm(T^+)$. Hence, 
\begin{align*}
 \frac{\log\left|\gT_r^\pm(B^H,T)-T\right|}{\log r} \geq \frac{\log(r^{1/H}+r^{o(1)})}{\log r} =\frac{1}{H}+o(1)
\end{align*}
and taking the limit inferior for $r\searrow0$ together with \eqref{eq:sumoflogs} yields
\[
 \liminf_{r\searrow0} \frac{\log \PTEp{\tau_{\gp(W)}(B_\infty(x,r))}{x}}{\log r} \geq \frac{2}{H},
\]
thus concluding the proof of Theorem~\ref{thm:BHdimw}.

In a more general setting, replacing $B^H$ by the graph of an arbitrary everywhere local $\ga$-H\"older regular map, we can by similar means obtain a concise statement for global upper walk dimension:
\begin{proposition}\label{prop:ugwdim}
 If $(X,d_X)$ is a metric space and $f:\IR\to X$ is everywhere locally $\ga$-H\"older regular, then 
 \[
  \dimuw(\gr(f),\gp(W))=\frac{2}{\ga},
 \]
 where we endow $\gr(f)$ with the restriction of the metric $d_\infty$ on $\IR\times X$, defined by $d_\infty((t,x),(t',x')=|t-t'|\wedge d_X(x,x')$.
\end{proposition}
\begin{proof}
 Since the inequality ``$\leq$'' is ensured by Lemma~\ref{lem:dimuw1}, and it remains to show ``$\geq$''.
 As before, consider $x=(T,f(T))$, and denote balls of radius $r$ around $x$ in the $d_\infty$-metric on $\gr(f)$ by $B_\infty(x,r)$. Analogously to \eqref{eq:Theta}, we define
 \begin{align*}
  \gT^+_r(f,T)&:=\inf\big\{t>T:(t,f(t))\notin B_\infty(x,r)\big\}\\
  \gT^-_r(f,T)&:=\sup\big\{t<T:(t,f(t))\notin B_\infty(x,r)\big\}.
 \end{align*}
 Hence, we obtain an equation analogous to \eqref{eq:sumoflogs} in the same way as in the proof of Theorem~\ref{thm:BHdimw}. Therefore, it remains to show that 
   \begin{equation}\label{eq:limtwo}
    \limsup_{r\searrow 0}\frac{\log\left|\gT^{\pm}_r(f,T)-T\right|}{\log r}\geq\frac{1}{\ga}\ \ \ \forall T\in\IR,
  \end{equation}
  where both inequalities can be shown independently. To this end, choose $1\geq\gb>\ga$ and $T\in\IR$ arbitrarily. Then, there exist sequences $T-r^-_n\nearrow T$ and $T+r^+_n\searrow T$ such that 
  \[
    \left|f(T)-f(T\pm r^{\pm}_n)\right|>(r^{\pm}_n)^\gb,
  \]
  because $f$ is nowhere locally $\gb$-H\"{o}lder continuous. Thus we have $\gT^{\pm}_{(r_n)^\gb}(f,T)\leq r_n$ which in turn implies
  \[
    \frac{\log \gT^{\pm}_{r_n^\gb}(f,T)}{\log r_n^\gb} 
    \geq \frac{1}{\gb},
  \]
  where $r_n:=r_n^{\pm}$ for brevity. Hence, 
  \[
    \limsup_{r\searrow 0}\frac{\log\left|\gT^{\pm}_r(f,T)-T\right|}{\log r}\geq\frac{1}{\gb}\ \ \ \forall T\in\IR
  \]
  and because the left-hand side does not depend on $\gb$, we can take the limit for $\gb\searrow H$ to obtain \eqref{eq:limtwo}. 
\end{proof}

\section{Further Questions}

As we saw, the Einstein relation is an invariant of metric measure spaces. Under H\"older continuous tranformations, its behaviour depends on the Hausdorff and walk dimensions, for both of which we have the same upper bound, but in general different behaviour. The previously presented results, such as the one about Brownian motion on the graph of a fractional Brownian motion, aside there remains a plethora of further questions to discuss such as (arranged in order of increasing speculativeness):

\textbf{Is there a general lemma providing lower estimates for the walk dimension?} This question is almost self-explanatory, and it aims at a statement that plays a similar role for the walk dimension as the mass-distribution principle does for the Hausdorff dimension.

\textbf{Is the Hausdorff dimension the ``right'' fractal dimension for the Einstein relation?} There are several alternative ways to define geometric dimensions for fractals, such as the packing dimension or the box-counting dimension. Of course, any reasonably well-behaved notion of dimension should be definable on a large class of metric spaces and should be invariant under isometries. Naturally, then the arguments used in the proof of Proposition~\ref{prop:mmiso} in the stricter setting of mm-isomorphisms show that the Einstein relation will still be an invariant of mm-spaces. 
  
When evaluating how ``good'' a fractal dimension for this purpose is, two questions should be asked: 
\begin{enumerate}
  \item Can this variant of the Einstein relation distinguish between spaces that are Lipschitz- yet not mm-isomorphic and if so, does it better than other variants?
  \item Are there general theorems for this variant that give explanations on why the Einstein relation should hold with constant $c$ on interesting classes of spaces?
\end{enumerate}
Of course, the latter questions are difficult to answer and not much is understood yet even for $\dimh$. 

A similar question is whether there exists a variation to the Einstein relation that can tell apart spaces that are Lipschitz- but not mm-isomorphic.

\textbf{Is it possible to extend the Einstein relation \eqref{eq:ER} to graphs in such a way that it is compatible with the discrete version from Section~7?}
We saw in Section~8 that the local walk dimension is better suited for bounded metric spaces. On the other hand, approximating spaces by a sequence of finite graphs as in the case of the Sierpi\'nski gasket is a useful tool to have. However, for graphs the limit $r\searrow0$ in the definition of the walk dimension does not make sense. 

One way to circumvent these problems with a unified approach might be to consider metric graphs $\scG$. Here, a metric graph is a disjoint collection of closed intervals $I_i$, where either $I_i=[a_i,b_i]$ or $I_i=[a_i,\infty)$ for $a_i,b_i\in\IR$, $i\in\scI$ an index set, together with an equivalence relation $\sim$ on the set of boundary points $\{a_i,b_i:i\in\scI\}$, where the boundary points are identified according to $\sim$. In other words, $\scG$ is the quotient space
\[
  \Big(\bigsqcup_{i\in\scI} I_i\Big)\Big/\!\sim.
\]
As stochastic processes on metric graphs have been investigated in recent years (cf. \cite{werner2016brownian}), it is a natural question to ask whether one can replace the approximation of the Brownian motion on $\SG$ by random walks on $G_n$ with an approximation by Brownian motions on $\scG_n$, where $\scG_n$ are the metric graphs with the metric structure coming from the embedding of $G_n$ in $\IR^2$. If this happens to be the case, one can furthermore ask if Definition~\ref{def:dimw}, applied to the approximating processes on $\scG_n$, yields an approximation of the walk dimension on $\SG$. 

\textbf{What are the topological properties of the Einstein relation?} This question aims at finding a general setting in which the Einstein relation on a given space can be approximated by Einstein relations on other spaces.

The class of isomorphism classes of (compact!) mm-spaces, 
$\mathsf{cMM}_{\leq1}\big/\!\cong$, can be endowed with different topologies, perhaps the most well-known way of doing this is via the Gromov-Hausdorff-Prohorov metric, defined in the following way: 
  
Let $(X,d_X,\mu_X)$ and $(Y,d_Y,\mu_Y)$ be compact mm-spaces. Denote by $d_H$ the usual Hausdorff distance between closed sets in a metric space (cf. section 1.1) and by $d_P$ the Prohorov distance between probability measures $\nu,\nu'$ on $Z$,
\[
  d_P(\nu,\nu'):=\inf\left\{\gep>0: \nu(B(A,\gep))>\nu'(A)-\gep\
    \text{ for any } A\in\scB(Z)\right\}.
\]  
Then,
\[
  d_{GHP}(Z;\gi_X,\gi_Y):=\inf_{\gi_X,\gi_Y}
     \Big(d_H\big(\gi_X(X),\gi_Y(Y)\big)
     +d_P\big((\gi_X)_*\mu_X,(\gi_Y)_*\mu_Y\big)\Big),
\]
where the infimum is taken over all metric spaces $(Z,d_Z)$ and all isometric embeddings $\gi_X:X\hookrightarrow Z$, 
$\gi_Y:Y\hookrightarrow Z$ of $X,Y$ into $Z$. This defines a pseudo-metric on $\mathsf{cMM}_{\leq1}\big/\!\cong$. 
  
Now take a subset $\mathsf{X}\ssq\mathsf{cMM}_{\leq1}$ and a mapping 
\[
  \mathsf{X}\ni(X,d_X,\mu_X) \mapsto (A_X,\scD(A_X))\in\IL(L^2(X,\mu_X))
\]
that assigns to each mm-space an operator that satisfies the conditions \ref{cond:A}. After taking the quotient, we are left with a map that sends each mm-isomorphism class to a linear operators on a representative of this class,
\[
  \left(\mathsf{X}\big/\!\cong\right) \ni [(X,d_X,\mu_X)]_{\cong}
  \mapsto (A_X,\scD(A_X)) \in \IL(L^2(X,\mu_X)),
\]
where the right-hand side is unique up to the transport of structure induced by mm-isomorphisms as discussed in section 9. In particular, we can now regard the constant in the Einstein relation as a function 
\begin{align}
  \scE:\mathsf{X}\big/\!\cong&\to \IR_{\geq0}\notag\\
  \left[(X,d_X,\mu_X)\right]_{\cong} &\mapsto 
    c=\frac{\dimh(X)}{\dims(X,A_X)\dimw(X,M^{A_X})}
\end{align}
In general, fixing a topology $\scT$ on $\mathsf{cMM}_{\leq1}\big/\!\cong$, a set $\mathsf{X}$ as above and an assignment 
$A_X$ (which should, in some sense, depend continuously on 
$[X]$), what can be said about the topological properties of 
$\scE$? Is it continuous w.r.t $\scT$? Are at least the preimages of single points closed?


\end{document}